\documentclass[english]{article}
\usepackage[T1]{fontenc}
\usepackage[latin9]{inputenc}
\usepackage{amsmath}
\usepackage{amsthm}
\usepackage{amssymb}

\makeatletter

\providecommand{\tabularnewline}{\\}

\theoremstyle{plain}
\newtheorem{thm}{\protect\theoremname}
  \theoremstyle{remark}
  \newtheorem{claim}[thm]{\protect\claimname}

\usepackage{fullpage}\usepackage{makeidx}\usepackage{amsfonts}\usepackage{latexsym}

\theoremstyle{definition}
\newtheorem{lemma}{Lemma}[section]
\newtheorem{theorem}[lemma]{Theorem}\newtheorem{proposition}[lemma]{Proposition}\newtheorem{definition}[lemma]{Definition}\newtheorem{remark}[lemma]{Remark}\usepackage{times}

\usepackage{babel}

\usepackage[blocks]{authblk}
\title{$6A$-Algebra and its representations}

\author{Chongying Dong\footnote{supported by China NSF grant 11871351}}
\affil{Department of Mathematics, University of
California, Santa Cruz, CA 95064 USA}
\author{Xiangyu Jiao\footnote{Supported by China NSF 11401213, 11571391, 11671138, Science and Technology Commission of Shanghai Municipality (STCSM) 18dz2271000, 16ZR1417800}}
\affil{Department of Mathematics, East China Normal University, Shanghai 200241, CHINA}
\author{Nina Yu\footnote{Supported by China NSF 11601452, Fundamental Research Funds for the Central Universities 20720170010, and  Research Fund for Fujian Young Faculty JAT170006}}
\affil{School of Mathematical Sciences, Xiamen University, Xiamen, Fujian 361005, CHINA}

\usepackage{babel}

\usepackage{babel}

\usepackage{babel}

\usepackage{babel}
\providecommand{\claimname}{Claim}
\providecommand{\theoremname}{Theorem}

\usepackage{babel}
\providecommand{\claimname}{Claim}
\providecommand{\theoremname}{Theorem}

\numberwithin{equation}{section}

\usepackage{babel}
\providecommand{\claimname}{Claim}
\providecommand{\theoremname}{Theorem}

\usepackage{babel}
\providecommand{\claimname}{Claim}
\providecommand{\theoremname}{Theorem}

\makeatother

\usepackage{babel}
  \providecommand{\claimname}{Claim}
\providecommand{\theoremname}{Theorem}

\begin{document}
\maketitle
\begin{abstract}
In this paper, we study the structure and representation of a $6A$-algebra
which is a vertex operator algebra generated by two Ising vectors
$e,f$ with inner product $\left\langle e,f\right\rangle =\frac{5}{2^{10}}.$
In particular, we prove the uniqueness of the vertex operator algebra
structure of this 6A-algebra, classify the irreducible modules, and
determine the fusion rules.
\end{abstract}

\section{{\normalsize{}{}Introduction}}

This paper is devoted to the study of the $6A$-algebra which is a
vertex operator algebra of the Moonshine type generated by two Ising
vectors whose inner product is $5/2^{10}.$

An \emph{Ising vector} in a vertex operator algebra $V$ is a Virasoro
vector which generates a subalgebra isomorphic to the Virasoro vertex
operator algebra $L(\frac{1}{2},0)$. The importance of the Ising
vectors was first noticed in \cite{DMZ} for understanding the structure
of the moonshine vertex operator algebra $V^{\natural}$ \cite{FLM}.
In fact, $V^{\natural}$ contains a conformal subalgebra $L(\frac{1}{2},0)^{\otimes48}.$
This led to the theory of framed vertex operator algebras \cite{M,DGH},
a new construction of the moonshine vertex operator algebra \cite{M0},
a proof of Frenkel-Lepowsky-Meurman's conjecture \cite{FLM} that
$V^{\natural}$ is holomorphic \cite{D} and two weaker versions of
Frenkel-Lepowsky-Meurma's uniqueness conjecture on $V^{\natural}$
\cite{DGL,LY}. Moreover for an Ising vector $e\in V$, one can define
the Miyamoto involution $\tau_{e}$ which is an automorphism of $V$.

The study of a vertex operator algebra generated by two Ising vectors
initiated in \cite{M1}. It was proved in \cite{M1} that each axis
of the Monster Griess algebra is essentially a half of an Ising vector
of $V^{\sharp}$ and $\tau_{e}$ is a $2A$-involution of the Monster
simple group $\mathbb{M}$. Thus there is a one-to-one correspondence
between $2A$-involutions of $\mathbb{M}$ and Ising vectors of $V^{\sharp}$.
It is shown in \cite{C} that the structure of the subalgebra generated
by two Ising vectors $e$ and $f$ in the algebra $V_{2}^{\sharp}$
depends on only the conjugacy class of $\tau_{e}\tau_{f}$, and the
inner product $\left\langle e,f\right\rangle $ is given by the following
table:
\begin{center}
\begin{tabular}{|c|c|c|c|c|c|c|c|c|c|}
\hline
$\left\langle \tau_{e}\tau_{f}\right\rangle ^{\mathbb{M}}$  & $1A$  & $2A$  & $3A$  & $4A$  & $5A$  & $6A$  & $3C$  & $4B$  & $2B$\tabularnewline
\hline
$\left\langle e,f\right\rangle $  & $1/4$  & $1/2^{5}$  & $13/2^{10}$  & $1/2^{7}$  & $3/2^{9}$  & $5/2^{10}$  & $1/2^{8}$  & $1/2^{8}$  & $0$\tabularnewline
\hline
\end{tabular}
\par\end{center}

Let $V$ be an arbitrary simple vertex operator algebra of the moonshine
type. It was proved in \cite{S} that the structure of the subalgebra
generated by two Ising vectors in the Griess algebra $V_{2}$ of $V$
is uniquely determined by the inner product of the two Ising vectors.
Moreover, the inner product of two Ising vectors again has 9 possibilities
as in the case of the Moonshine vertex operator algebra. Certain vertex
operator subalgebras $\mathcal{U}_{nX}$ of the lattice vertex operator
algebra $V_{\sqrt{2}E_{8}}$ corresponding to the type $nX$ of $\left\langle \tau_{e}\tau_{f}\right\rangle ^{\mathbb{M}}$
were constructed in \cite{LYY1}. It was shown that in each of the
nine cases $\mathcal{U}_{nX}$ always contains conformal vectors $\hat{e}$
and $\hat{f}$ of central charge $1/2$ such that the inner product
$\langle\hat{e},\hat{f}\rangle$ is exactly those given in the table.
The structure and representation of these coset subalgebras $\mathcal{U}_{nX}$
are studied in \cite{LYY2} and it was shown that they are all generated
by two conformal vectors of central charge $1/2.$ It is also shown
that the product of two Miyamoto involutions is in the desired conjugacy
class of the Monster simple group if a coset subalgebra $\mathcal{U}_{nX}$
is actually contained in the Moonshine vertex operator algebra $V^{\sharp}.$
The existence of $\mathcal{U}_{nX}$ inside the Moonshine vertex operator
algebra $V^{\sharp}$ for the cases $nX=1A,2A,2B$ and $4A$ is also
established. Furthermore, the cases for $3A$, $5A$ and $3C$ are
discussed in \cite{LYY2}.

But the structure and representation of $\mathcal{U}_{6A}$ has not
been understood well. It turns out that this $6A$-algebra $\mathcal{U}_{6A}=U^{1}\oplus U^{2}\oplus U^{3}$
is an extension of a rational, $C_{2}$-cofinite vertex operator algebra
$U^{1}$ by two irreducible $U^{1}$-modules $U^{2}$ and $U^{3}$
which are not simple current modules. The first goal is to establish
the uniqueness of the $6A$-algebra. The main idea is to use relevant
braiding matrices. The second goal is to classify irreducible modules
for $\mathcal{U}_{6A}$, we first construct 14 irreducible $\mathcal{U}_{6A}$-modules
and then prove the sum of squares of quantum dimensions of these irreducible
modules is exactly the global dimension of $\mathcal{U}_{6A}$. This
implies that $\mathcal{U}_{6A}$ has exactly 14 irreducible modules.
Last, we determine the fusion rules of these modules.

For simplicity we denote $\mathcal{U}_{6A}$ by $\mathcal{U}$.

The paper is organized as follows. In Section 2, we review some basic
notions and some well known results in the vertex operator algebra
theory. In Section 3, we study the structure of the $6A$-algebra
and prove the uniqueness of the vertex operator algebra structure
on $\mathcal{U}$. In section 4, we classify the irreducible modules
for $\mathcal{U}$. In section 5, we determine the fusion rules among
irreducible $\mathcal{U}$-modules.

\section{Preliminary}

In this section, we review the basics on vertex operators algebras,
the theory of quantum dimensions from \cite{DJX}, the coset realization
of the discrete series of the unitary representations for the Virasoro
algebra \cite{GKO} and the braiding matrices for certain Virasoro
vertex operator algebras \cite{FFK}.

\subsection{Basics}

Let $V=(V,\ Y,\ \mathbf{1},\ \omega)$ be a vertex operator algebra.
Let $Y(v,\ z)=\sum_{n\in\mathbb{Z}}v_{n}z^{-n-1}$ denote the vertex
operator of $V$ for $v\in V$, where $v_{n}\in\mbox{End}(V)$. We
first recall some basic notions from \cite{FLM,Z,DLM1,DLM3}.

\begin{definition} A vector $v\in V_{2}$ is called a \emph{Virasoro
vector with central charge $c_{v}$ }if it satisfies\emph{ $v_{1}v=2v$
}and $v_{3}v=\frac{c_{v}}{2}\mathbf{1}$. Then the operators $L_{n}^{v}:=v_{n+1},\ n\in\mathbb{Z}$,
satisfy the Virasoro commutation relation
\[
\left[L_{m}^{v},\ L_{n}^{v}\right]=\left(m-n\right)L_{m+n}^{v}+\delta_{m+n,\ 0}\frac{m^{3}-m}{12}c_{v}
\]
for $m,\ n\in\mathbb{Z}.$ A Virasoro vector $v\in V_{2}$ with central
charge $1/2$ is called an \emph{Ising vector }if $v$ generates the
Virasoro vertex operator algebra $L(1/2,\ 0)$.

\end{definition}

\begin{definition} An \emph{automorphism} $g$ of a vertex operator
algebra $V$ is a linear isomorphism of $V$ satisfying $g\left(\omega\right)=\omega$
and $gY\left(v,z\right)g^{-1}=Y\left(gv,z\right)$ for any $v\in V$.
We denote by $\mbox{Aut}\left(V\right)$ the group of all automorphisms
of $V$. \end{definition}

For a subgroup $G\le\mbox{Aut}\left(V\right)$ the fixed point set
$V^{G}=\left\{ v\in V|g\left(v\right)=v,\forall g\in G\right\} $
has a vertex operator algebra structure.

Let $g$ be an automorphism of a vertex operator algebra $V$ of order
$T$. Denote the decomposition of $V$ into eigenspaces of $g$ as:

\[
V=\oplus_{r\in\mathbb{Z}/T\text{\ensuremath{\mathbb{Z}}}}V^{r}
\]
where $V^{r}=\left\{ v\in V|gv=e^{2\pi ir/T}v\right\} $.

\begin{definition} A \emph{weak $g$-twisted $V$-module} $M$ is
a vector space with a linear map
\begin{align*}
Y_{M}: & V\to\left(\text{End}M\right)\{z\}\\
 & v\mapsto Y_{M}\left(v,z\right)=\sum_{n\in\mathbb{Q}}v_{n}z^{-n-1}\ \left(v_{n}\in\mbox{End}M\right)
\end{align*}
which satisfies the following: for all $0\le r\le T-1$, $u\in V^{r}$,
$v\in V$, $w\in M$,

\[
Y_{M}\left(u,z\right)=\sum_{n\in-\frac{r}{T}+\mathbb{Z}}u_{n}z^{-n-1},
\]

\[
u_{l}w=0\ for\ l\gg0,
\]

\[
Y_{M}\left(\mathbf{1},z\right)=Id_{M},
\]

\[
z_{0}^{-1}\text{\ensuremath{\delta}}\left(\frac{z_{1}-z_{2}}{z_{0}}\right)Y_{M}\left(u,z_{1}\right)Y_{M}\left(v,z_{2}\right)-z_{0}^{-1}\delta\left(\frac{z_{2}-z_{1}}{-z_{0}}\right)Y_{M}\left(v,z_{2}\right)Y_{M}\left(u,z_{1}\right)
\]

\begin{equation}
=z_{2}^{-1}\left(\frac{z_{1}-z_{0}}{z_{2}}\right)^{-r/T}\delta\left(\frac{z_{1}-z_{0}}{z_{2}}\right)Y_{M}\left(Y\left(u,z_{0}\right)v,z_{2}\right),\label{Jacobi for twisted V-module}
\end{equation}
where $\delta\left(z\right)=\sum_{n\in\mathbb{Z}}z^{n}$.

\end{definition}

\begin{definition}

A $g$-\emph{twisted $V$-module} is a weak $g$-twisted $V$-module
$M$ which carries a $\mathbb{C}$-grading induced by the spectrum
of $L(0)$ where $L(0)$ is the component operator of $Y(\omega,z)=\sum_{n\in\mathbb{Z}}L(n)z^{-n-2}.$
That is, we have $M=\bigoplus_{\lambda\in\mathbb{C}}M_{\lambda},$
where $M_{\lambda}=\left\{ w\in M|L(0)w=\lambda w\right\} $. Moreover,
$\dim M_{\lambda}$ is finite and for fixed $\lambda,$ $M_{\frac{n}{T}+\lambda}=0$
for all small enough integers $n.$

\end{definition}

\begin{definition}An \emph{admissible $g$-twisted $V$-module} $M=\oplus_{n\in\frac{1}{T}\mathbb{Z}_{+}}M\left(n\right)$
is a $\frac{1}{T}\mathbb{Z}_{+}$-graded weak $g$-twisted module
such that $u_{m}M\left(n\right)\subset M\left(\mbox{wt}u-m-1+n\right)$
for homogeneous $u\in V$ and $m,n\in\frac{1}{T}\mathbb{Z}.$

\end{definition}

If $g=Id_{V}$ we have the notions of weak, ordinary and admissible
$V$-modules \cite{DLM2}.

\begin{definition}A vertex operator algebra $V$ is called \emph{$g$-rational}
if the admissible $g$-twisted module category is semisimple. $V$
is called \emph{rational} if $V$ is $1$-rational. \end{definition}

It was proved in \cite{DLM2} that if $V$ is rational then there
are only finitely irreducible admissible $V$-modules up to isomorphism
and each irreducible admissible $V$-module is ordinary. Let $M^{0},M^{1},$
$\cdots,M^{d}$ be all the irreducible modules up to isomorphism with
$M^{0}=V$. Then there exists $h_{i}\in\mathbb{C}$ for $i=0,\cdots,d$
such that
\[
M^{i}=\oplus_{n=0}^{\infty}M_{h_{i}+n}^{i}
\]
where $M_{h_{i}}^{i}\not=0$ and $L\left(0\right)|_{M_{h_{i}+n}^{i}}=h_{i}+n$,
$\forall n\in\mathbb{Z}_{+}$. $h_{i}$ is called the \emph{conformal
weight} of $M^{i}$. We denote $M^{i}\left(n\right)=M_{h_{i}+n}^{i}.$

Let $M=\oplus_{\lambda_{\in\mathbb{C}}}M_{\lambda}$ be a $V$-module.
The restricted dual of $M$ is defined by $M'=\oplus_{\lambda_{\in\mathbb{C}}}M_{\lambda}^{\ast}$
where $M_{\lambda}^{\ast}=\text{Hom}_{\mathbb{C}}\left(M_{\lambda},\mathbb{C}\right).$
It was proved in \cite{FHL} that $M'=\left(M',Y_{M'}\right)$ is
naturally a $V$-module such that
\[
\left\langle Y_{M'}\left(v,z\right)f,u\right\rangle =\left\langle f,Y_{M}\left(e^{zL\left(1\right)}\left(-z^{-2}\right)^{L\left(0\right)}v,z^{-1}\right)u\right\rangle ,
\]
for $v\in V,$ $f\in M'$ and $u\in M$, and $\left(M'\right)'\cong M$.
Moreover, if $M$ is irreducible, so is $M'$. A $V$-module $M$
is said to be \emph{self dual} if $M\cong M'$.

\begin{definition} A vertex operator algebra $V$ is said to be \emph{$C_{2}$-cofinite}
if $V/C_{2}(V)$ is finite dimensional, where $C_{2}(V)=\langle v_{-2}u|v,u\in V\rangle.$
\end{definition}

\begin{definition}A vertex operator algebra $V=\oplus_{n\in\mathbb{Z}}V_{n}$
is said to be of \emph{CFT type} if $V_{n}=0$ for negative $n$ and
$V_{0}=\mathbb{C}1$. \end{definition}

\begin{definition} Let $\left(V,Y\right)$ be a vertex operator algebra
and let $\left(M^{i},Y^{i}\right),\ \left(M^{j},Y^{j}\right)$ and
$\left(M^{k},Y^{k}\right)$ be three $V$-modules. An \emph{intertwining
operator of type $\left(\begin{array}{c}
M^{k}\\
M^{i}\ M^{j}
\end{array}\right)$} is a linear map
\begin{gather*}
\mathcal{Y}\left(\cdot,z\right):\ M^{i}\to\text{\ensuremath{\mbox{Hom}\left(M^{j},\ M^{k}\right)\left\{ z\right\} }}\\
u\mapsto\mathcal{Y}\left(u,z\right)=\sum_{n\in\mathbb{Q}}u_{n}z^{-n-1}
\end{gather*}
satisfying:

(1) for any $u\in M^{i}$ and $v\in M^{j}$, $u_{n}v=0$ for $n$
sufficiently large;

(2) $\mathcal{Y}(L_{-1}v,\ z)=\left(\frac{d}{dz}\right)\mathcal{Y}\left(v,z\right)$
for $v\in M^{i}$;

(3) (Jacobi Identity) for any $u\in V,\ v\in M^{i}$,
\begin{alignat*}{1}
 & z_{0}^{-1}\delta\left(\frac{z_{1}-z_{2}}{z_{0}}\right)Y^{k}\left(u,z_{1}\right)\mathcal{Y}\left(v,z_{2}\right)-z_{0}^{-1}\delta\left(\frac{-z_{2}+z_{1}}{z_{0}}\right)\mathcal{Y}\left(v,z_{2}\right)Y^{j}\left(u,z_{1}\right)\\
 & =z_{2}^{-1}\left(\frac{z_{1}-z_{0}}{z_{2}}\right)\mathcal{Y}\left(Y^{i}\left(u,z_{0}\right)v,z_{2}\right).
\end{alignat*}
The space of all intertwining operators of type $\left(\begin{array}{c}
M^{k}\\
M^{i}\ M^{j}
\end{array}\right)$ is denoted $I_{V}\left(\begin{array}{c}
M^{k}\\
M^{i}\ M^{j}
\end{array}\right)$. Without confusion, we also denote it by $I_{i,j}^{k}.$ Let $N_{i,\ j}^{k}=\dim I_{i,j}^{k}$.
These integers $N_{i,j}^{k}$ are called the \emph{fusion rules}.
\end{definition}

The following proposition was proved in \cite{ADL}:

\begin{proposition} \label{restriction of fusion rules} Let $V$
be a vertex operator algebra and let $W^{1}$, $W^{2}$, $W^{3}$
be $V$-modules among which $W^{1}$ and $W^{2}$ are irreducible.
Suppose that $V_{0}$ is a vertex operator subalgebra of $V$ (with
the same Virasoro element) and that $N^{1}$ and $N^{2}$ are irreducible
$V_{0}$-modules of $W^{1}$ and $W^{2}$, respectively. Then the
restriction map from $I_{V}\left(_{W^{1}\ W^{2}}^{\ \ W^{3}}\right)$
to $I_{V_{0}}\left(_{N^{1}\ N^{2}}^{\ \ W^{3}}\right)$ is injective.
In particular,
\[
\dim I_{V}\left(_{W^{1\ }W^{2}}^{\ \ W^{3}}\right)\le\dim I_{V_{0}}\left(_{N^{1}\ N^{2}}^{\ \ W^{3}}\right).
\]
\end{proposition}

Let $V^{1}$ and $V^{2}$ be vertex operator algebras. Let $M^{i}$
, $i=1,2,3$, be $V^{1}$-modules, and $N^{i}$, $i=1,2,3$, be $V^{2}$-modules.
Then $M^{i}\otimes N^{i}$, $i=1,2,3$, are $V^{1}\otimes V^{2}$-modules
by \cite{FHL}. The following property was given in \cite{ADL}:

\begin{proposition} \label{fusion of tensor product}If $N_{M^{1},M^{2}}^{M^{3}}<\infty$
or $N_{N^{1},N^{2}}^{N^{3}}<\infty,$ then
\[
N_{M^{1}\otimes N^{1},M^{2}\otimes N^{2}}^{M^{3}\otimes N^{3}}=N_{M^{1},M^{2}}^{M^{3}}N_{N^{1},N^{2}}^{N^{3}}.
\]
\end{proposition}

Let $M^{1}$ and $M^{2}$ be $V$-modules. A fusion product for the
ordered pair $\left(M^{1},M^{2}\right)$ is a pair $\left(M,F\left(\cdot,z\right)\right)$
which consists of a $V$-module $M$ and an intertwining operator
$F\left(\cdot,z\right)$ of type $\left(_{M^{1},M^{2}}^{M}\right)$
such that the following universal property holds: For any $V$-module
$W$ and any intertwining operator $I\left(\cdot,z\right)$ of type
$\left(_{M^{1},M^{2}}^{W}\right)$, there exists a unique $V$-homomorphism
$\phi$ from $M$ to $W$ such that $I\left(\cdot,z\right)=\phi\circ F\left(\cdot,z\right).$
It is clear from the definition that if a tensor product of $M^{1}$
and $M^{2}$ exists, it is unique up to isomorphism. In this case,
we denote the \emph{fusion product} by $M^{1}\boxtimes_{V}M^{2}.$

The basic result is that the fusion product exists if $V$ is rational.
Let $M,N$ be irreducible $V$ -modules, we shall often consider the
fusion product
\[
M\boxtimes_{V}N=\sum_{W}N_{M,\ N}^{W}W
\]
where $W$ runs over the set of equivalence classes of irreducible
$V$-modules.

\begin{definition} Let $V$ be a simple vertex operator algebra.
A simple $V$-module $M$ is called a \emph{simple current} if for
any irreducible $V$-module $W$, $M\boxtimes_{V}W$ exists and is
also a simple $V$-module. \end{definition}

The following proposition is from \cite{FHL}:

\begin{proposition} \label{extension property} Let $V$ be a vertex
operator algebra and $V'$ be its restricted dual. For $u,v,w\in V$
and $t\in V'$, we have the following equality of rational functions

\begin{equation}
\iota_{12}^{-1}\left\langle t,Y\left(u,z_{1}\right)Y\left(v,z_{2}\right)w\right\rangle =\iota_{21}^{-1}\left\langle t,Y\left(v,z_{2}\right)Y\left(u,z_{1}\right)w\right\rangle \label{communitivity}
\end{equation}

\begin{equation}
\iota_{12}^{-1}\left\langle t,Y\left(u,z_{1}\right)Y\left(v,z_{2}\right)w\right\rangle =\left(\iota_{120}^{-1}\left\langle t,Y\left(Y\left(u,z_{0}\right)v,z_{2}\right)w\right\rangle \right)|_{z_{0}=z_{1}-z_{2}}\label{associavity}
\end{equation}
where $\iota_{12}^{-1}f\left(z_{1},z_{2}\right)$ denotes the formal
power expansion of an analytic function $f\left(z_{1},z_{2}\right)$
in the domain $\left|z_{1}\right|>\left|z_{2}\right|$ . \end{proposition}

The following result about bilinear form on $V$ is from \cite{Li2}:

\begin{theorem}\label{bilinear form}The space of invariant bilinear
forms on $V$ is isomorphic to the space
\[
\left(V_{0}/L\left(1\right)V_{1}\right)^{*}=\mbox{Hom}_{\mathbb{C}}\left(V_{0}/L\left(1\right)V_{1},\mathbb{C}\right).
\]
\end{theorem}

\subsection{Quantum Galois Theory}

Now we recall quantum Galois theory and quantum dimensions from \cite{DM}
and \cite{DJX}. Let $V$ be a simple vertex operator algebra and
$G$ a finite and faithful group of automorphisms of $V$. Let $\text{Irr}\left(G\right)$
be the set of simple characters $\chi$ of $G$. As $\mathbb{C}G$-module,
each homogeneous space $V_{n}$ of $V$ is finite dimensional, and
$V$ can be decomposed into a direct sum of graded subspaces
\[
V=\oplus_{\chi\in\text{Irr}\left(G\right)}V^{\chi},
\]
where $V^{\chi}$ is the subspace of $V$ on which $G$ acts according
to the character $\chi$. The following theorem is from \cite{DM}.

\begin{theorem} \label{classical galois theory}Suppose that $V$
is a simple vertex operator algebra and that $G$ is a finite and
faithful solvable group of automorphisms of $V$. Then the following
hold:

(i) Each $V^{\chi}$ is nonzero;

(ii) For $\chi\in\text{Irr}\left(G\right)$, each $V^{\chi}$ is a
simple module for the $G$-graded vertex operator algebra $\mathbb{C}G\otimes V^{G}$
of the form

\[
V^{\chi}=M_{\chi}\otimes V_{\chi},
\]
where $M_{\chi}$ is the simple $\mathbb{C}G$-module affording $\chi$
and where $V_{\chi}$ is a simple $V^{G}$-module.

(iii) The map $M_{\chi}\mapsto V_{\chi}$ is a bijection from the
set of simple $\mathbb{C}G$-modules to the set of (inequivalent)
simple $V^{G}$-modules which are contained in $V$.

\end{theorem}

Now we recall the notion of quantum dimension from \cite{DJX}. Let
$V$ be a vertex operator algebra of CFT type and $M$ a $V$-module,
the formal character of $M=\oplus_{n\in\mathbb{Z}_{+}}M_{\lambda+n}$
is defined to be
\[
Ch_{q}M=\mbox{tr}_{q}M=\text{tr}q^{L(0)-c/24}=q^{\lambda-c/24}\sum_{n\in\mathbb{Z}_{+}}(\dim M_{\lambda+n})q^{n}
\]
where $\lambda$ is the conformal weight of $M$. The quantum dimension
of $M$ over $V$ is defined as:
\[
q\dim_{V}M=\lim_{q\to1^{-}}\frac{Ch_{q}M}{Ch_{q}V}.
\]
The following result is from Theorem 6.3 in \cite{DJX}:

\begin{theorem} \label{quantum dimension and orbifold module}Let
$V$ be a rational and $C_{2}$-cofinite simple vertex operator algebra.
Assume $V$ is $g$-rational and the conformal weight of any irreducible
$g$-twisted $V$-module is positive except for $V$ itself for all
$g\in G$. Then
\[
q\dim_{V^{G}}V_{\chi}=\dim W_{\chi}.
\]
\end{theorem}

For convenience, from now on, we say a vertex operator algebra $V$
is ``good'' if it satisfies the following conditions: $V$ is a
rational and $C_{2}$-cofinite simple vertex operator algebra of CFT
type with $V\cong V'$. Let $M^{0},\ M^{1},\ \cdots,\ M^{d}$ be all
the inequivalent irreducible $V$-modules with $M^{0}\cong V$ . The
corresponding conformal weights $\lambda_{i}$ satisfy $\lambda_{i}>0$
for $0<i\le d$.

The following properties of quantum dimensions are from \cite{DJX}
:

\begin{proposition} \label{product property of quantum dimension }\label{quantum dimension and simple current}
Let $V$ be a ``good'' vertex operator algebra. Then

(i) $q\dim_{V}\left(M^{i}\boxtimes M^{j}\right)=q\dim_{V}M^{i}\cdot q\dim_{V}M^{j},$
$\forall i,j$.

(ii) A $V$-module $M^{i}$ is a simple current if and only if $q\dim_{V}M^{i}=1$.

(iii) $q\dim_{V}M^{i}\in\left\{ 2\cos\left(\pi/n\right)|n\ge3\right\} \cup\left\{ a|2\le a<\infty,a\ is\ algebraic\right\} .$

\end{proposition}

\begin{definition} Let $V$ be a vertex operator algebra with finitely
many inequivalent irreducible modules $M^{0},\cdots,M^{d}$. The \emph{global
dimension} of $V$ is defined as
\[
\text{glob}\left(V\right)=\sum_{i=0}^{d}\left(q\dim_{V}M^{i}\right)^{2}.
\]
\end{definition}

\begin{remark} \label{product property of global dimension}Let $U$
and $V$ be ``good'' vertex operator algebras, $M$ be a $U$-module
and $N$ be a $V$-module. Then Lemma 2.10 of \cite{ADJR} gives

\[
q\dim_{U\otimes V}M\otimes N=q\dim_{U}M\cdot q\dim_{U}N,
\]
\[
\text{glob}\left(U\otimes V\right)=\text{glob}\left(U\right)\cdot\text{glob}\left(V\right).
\]
\end{remark}

Let $V$ be a vertex operator algebra, recall that a simple vertex
operator algebra containing $V$ is called an \emph{extension} $U$
of $V$. Now we have the following theorem \cite{ABD,HKL,ADJR}:

\begin{theorem} \label{rationality of extesnion of VOA } \label{global dim of VOA and subVOA}
Let $V$ be a ``good'' vertex operator algebra. Let $U$ be a simple
vertex operator algebra which is an extension of $V$. Then $U$ is
also ``good'' and
\[
\text{glob}\left(V\right)=\text{glob}\left(U\right)\cdot\left(q\dim_{V}\left(U\right)\right)^{2}.
\]

\end{theorem}

\subsection{\label{subsec:The-unitary-series}The unitary series of the Virasoro
VOAs}

Now we recall notations about unitary minimal models of Virasoro algebra
from \cite{FFK}. The models are parameterized by a complex number
$\alpha_{-}^{2}$, related to the central charge of the Virasoro algebra
by $c=13-6\alpha_{-}^{2}-6\alpha_{-}^{-2}$ where $\alpha_{-}^{2}=\frac{p}{p'}$
and $\left|p-p'\right|=1.$ Without loss of generality, we write $p'=p+1$
and denote $c_{p}=1-\frac{6}{p\left(p+1\right)}$ with $p=2,3,4,\cdots$.
The label $I$ stands for a pair $\left(i',i\right)$ of positive
integers and the corresponding highest weight is
\begin{equation}
h_{I}=h_{\left(i'\ i\right)}^{\left(p\right)}=\frac{1}{4}\left(i'^{2}-1\right)\alpha_{-}^{2}-\frac{1}{2}\left(i'i-1\right)+\frac{1}{4}\left(i^{2}-1\right)\alpha_{-}^{-2}=\frac{\left(pi'-\left(p+1\right)i\right)^{2}-1}{4p\left(p+1\right)}.\label{minimal model notation}
\end{equation}
for $1\le i'\le p,$ $1\le i\le p-1.$ We denote such unitary minimal
models of Virasoro algebra by $L\left(c_{p},h_{\left(i',i\right)}^{\left(p\right)}\right).$

\begin{remark} \label{pairs corr. to U's weights for unitary model}Use
the above notation, we see that the central charge of the model $L\left(\frac{25}{28},0\right)$
corresponds to the parameter $\alpha_{-}^{2}=\frac{7}{8}$ with $p=7,$
$p'=8$. The highest weights for irreducible $L\left(\frac{25}{28},0\right)$-modules
are

\begin{equation}
\left\{ 0,\frac{5}{32},\frac{3}{4},\frac{57}{32},\frac{13}{4},\frac{165}{32},\frac{15}{2},\frac{5}{14},\frac{3}{224},\frac{3}{28},\frac{143}{224},\frac{45}{28},\frac{675}{224},\frac{34}{7},\frac{9}{7},\frac{99}{224},\frac{1}{28},\frac{15}{224},\frac{25}{28},\frac{323}{224},\frac{39}{14}\right\} .\label{highest weights for 25/28}
\end{equation}
In particular, the pairs $\left(1,1\right)$, $\left(1,5\right)$
and $\left(1,3\right)$ correspond to the highest weights 0, $\frac{34}{7}$
and $\frac{9}{7}$ respectively. \end{remark}

Also note that the fusion rules for irreducible $L\left(c_{p},0\right)$-modules
are as follows \cite{W}:

\begin{definition} An ordered triple of pairs of integers $\left(\left(i',i\right),\left(j',j\right),\left(k',k\right)\right)$
is called \emph{admissible }if $1\le i',j',k'\le p+1,1\le i,j,k\le p$,
$i'+j'+k'<2\left(p+1\right),$ $i+j+k<2p$, $i'<j'+k'$, $j'<i'+k'$,
$k'<i'+j'$, $i<j+k$, $j<i+k$, $k<i+j$, and the sums $i'+j'+k'$,
$i+j+k$ are odd. \end{definition}

\begin{proposition} \label{fusion rules of virasoro modules}The
fusion rules between $L\left(c_{p},0\right)$-modules $L\left(c_{p},h_{\left(i',i\right)}^{\left(p\right)}\right),$
$L\left(c_{p},h_{\left(j',j\right)}^{\left(p\right)}\right)$ are
\[
L\left(c_{p},h_{\left(i',i\right)}^{\left(p\right)}\right)\boxtimes L\left(c_{p},h_{\left(j',j\right)}^{\left(p\right)}\right)=\sum_{\left(k',k\right)}N_{\left(i',i\right),\left(j',j\right)}^{\left(k',k\right)}L\left(c_{p},h_{\left(k',k\right)}^{\left(p\right)}\right),
\]
where $N_{\left(i',i\right),\left(j',j\right)}^{\left(k',k\right)}$
is $1$ iff $\left(\left(i',i\right),\left(j',j\right),\left(k',k\right)\right)$
is an admissible triple of pairs and $0$ otherwise. \end{proposition}

\subsection{\label{subsec:Braiding-matrix}Braiding matrices}

Now we recall four point functions. Let $V$ be a rational and $C_{2}$-cofinite
vertex operator algebra of CFT type and $V\cong V'$. Let $M^{a_{1}},M^{a_{2}},M^{a_{3}},M^{a_{4}}$
be four irreducible $V$-modules. By Lemma 4.1 in \cite{H2} we know
that for $u_{a_{i}}\in M^{a_{i}},$

\[
\left\langle u_{a_{4}'},\mathcal{Y}_{a_{1},a_{5}}^{a_{4}}\left(u_{a_{1}},z_{1}\right)\mathcal{Y}_{a_{2},a_{3}}^{a_{5}}\left(u_{a_{2}},z_{2}\right)u_{a_{3}}\right\rangle ,
\]

\[
\left\langle u_{a_{4}'},\mathcal{Y}_{a_{2},a_{6}}^{a_{4}}\left(u_{a_{1}},z_{2}\right)\mathcal{Y}_{a_{1},a_{3}}^{a_{6}}\left(u_{a_{1}},z_{1}\right)u_{a_{3}}\right\rangle
\]
are analytic on $\left|z_{1}\right|>\left|z_{2}\right|>0$ and $\left|z_{2}\right|>\left|z_{1}\right|>0$
respectively, and can both be analytically extended to multi-valued
analytic functions on
\[
R=\left\{ \left(z_{1},z_{2}\right)\in\mathbb{C}^{2}|z_{1},z_{2}\not=0,z_{1}\not=z_{2}\right\} .
\]
One can lift the multi-valued functions on $R$ to single-valued functions
on the universal covering $\tilde{R}$ to $R$ as in \cite{H3}. We
use

\[
E\left\langle u_{a_{4}'},\mathcal{Y}_{a_{1},a_{5}}^{a_{4}}\left(u_{a_{1}},z_{1}\right)\mathcal{Y}_{a_{2},a_{3}}^{a_{5}}\left(u_{a_{2}},z_{2}\right)u_{a_{3}}\right\rangle
\]
and
\[
E\left\langle u_{a_{4}'},\mathcal{Y}_{a_{2},a_{6}}^{a_{4}}\left(u_{a_{1}},z_{2}\right)\mathcal{Y}_{a_{1},a_{3}}^{a_{6}}\left(u_{a_{1}},z_{1}\right)u_{a_{3}}\right\rangle
\]
to denote those analytic functions.

Let $\left\{ \mathcal{Y}_{a,b;i}^{c}|i=1,\cdots,N_{a,b}^{c}\right\} $
be a basis of $I_{a,b}^{c}$. From \cite{H3},

\[
\left\{ E\left\langle u_{a_{4}'},\mathcal{Y}_{a_{1},a_{5};i}^{a_{4}}\left(u_{a_{1}},z_{1}\right)\mathcal{Y}_{a_{2},a_{3};j}^{a_{5}}\left(u_{a_{2}},z_{2}\right)u_{a_{3}}\right\rangle |i=1,\cdots,N_{a_{1},a_{5}}^{a_{4}},j=1,\cdots,N_{a_{2},a_{3}}^{a_{5}},\forall a_{5}\right\}
\]
is a linearly independent set. Fix a basis of intertwining operators.
It was proved in \cite{KZ,TK} that
\begin{alignat*}{1}
 & span\left\{ E\left\langle u_{a_{4}'},\mathcal{Y}_{a_{3},\mu;i}^{a_{4}}\left(u_{a_{3}},z_{1}\right)\mathcal{Y}_{a_{2},a_{1};j}^{\mu}\left(u_{a_{2}},z_{2}\right)u_{a_{1}}\right\rangle |i,j,\mu\right\} \\
 & =span\left\{ E\left\langle u_{a_{4}'},\mathcal{Y}_{a_{2},\gamma;k}^{a_{4}}\left(u_{a_{2}},z_{2}\right)\mathcal{Y}_{a_{3},a_{1};l}^{\gamma}\left(u_{a_{3}},z_{1}\right)u_{a_{1}}\right\rangle |k,l,\gamma\right\} ,
\end{alignat*}
where $u_{a_{i}}\in M^{a_{i}}$. Then there exists $\left(B_{a_{4,}a_{1}}^{a_{3},a_{2}}\right)_{\mu,\gamma}^{i,j;k,l}\in\mathbb{C}$
such that
\begin{alignat}{1}
 & E\left\langle u_{a_{4}'},\mathcal{Y}_{a_{3},\mu;i}^{a_{4}}\left(u_{a_{3}},z_{1}\right)\mathcal{Y}_{a_{2},a_{1};j}^{\mu}\left(u_{a_{2}},z_{2}\right)u_{a_{1}}\right\rangle \nonumber \\
 & =\sum_{k,l,\gamma}\left(B_{_{a_{4},a_{1}}}^{a_{3,}a_{2}}\right)_{\mu,\gamma}^{i,j;k,l}E\left\langle u_{a_{4}'},\mathcal{Y}_{a_{2},\gamma;k}^{a_{4}}\left(u_{a_{2}},z_{1}\right)\mathcal{Y}_{a_{3},a_{1};l}^{\gamma}\left(u_{a_{3}},z_{2}\right)u_{a_{1}}\right\rangle \label{braiding matrix property}
\end{alignat}
(see \cite{H1,H2}). $B_{_{a_{4},a_{1}}}^{a_{3,}a_{2}}$ is called
the\emph{ braiding matrix}.

Let $b_{1},b_{2},b_{3},b_{4}$ be four irreducible $L\left(\frac{25}{28},0\right)$-modules.
Fix a basis $\left\{ \mathcal{\overline{Y}}_{a,b;i}^{c}|i=1,\cdots,N_{a,b}^{c}\right\} $
of intertwining operators of $I_{L\left(\frac{25}{28},0\right)}\left(_{L\left(\frac{25}{28},a\right),L\left(\frac{25}{28},b\right)}^{L\left(\frac{25}{28},c\right)}\right)$
with $N_{a,b}^{c}=\dim I_{L\left(\frac{25}{28},0\right)}\left(_{L\left(\frac{25}{28},a\right),L\left(\frac{25}{28},b\right)}^{L\left(\frac{25}{28},c\right)}\right)$
as in \cite{FFK}. Then there exists a matrix $\left(\tilde{B}_{_{b_{4},b_{1}}}^{b_{3},b_{2}}\right)_{\mu,\gamma}^{i,j;k,l}\in\mathbb{C}$
such that
\begin{alignat}{1}
 & E\left\langle u_{b_{4}'},\mathcal{\overline{Y}}_{b_{3},\mu;i}^{b_{4}}\left(u_{b_{3}},z_{1}\right)\mathcal{\overline{Y}}_{b_{2},b_{1};j}^{\mu}\left(u_{b_{2}},z_{2}\right)u_{b_{1}}\right\rangle \nonumber \\
 & =\left(\tilde{B}_{_{b_{4},b_{1}}}^{b_{3,}b_{2}}\right)_{\mu,\gamma}^{i,j;k,l}E\left\langle u_{b_{4}'},\mathcal{\overline{Y}}_{b_{2},\gamma;k}^{b_{4}}\left(u_{b_{2}},z_{2}\right)\mathcal{\overline{Y}}_{b_{3},b_{1};l}^{\gamma}\left(u_{b_{3}},z_{1}\right)u_{b_{1}}\right\rangle \label{tilde=00007BB=00007D}
\end{alignat}
by (\ref{braiding matrix property}).Now we recall some formulas about
minimal models of Virasoro vertex operator algebra given in \cite{FFK}.
We will use these formulas to prove some properties of braiding matrices,
which will be needed in the proof of uniqueness of the structure of
the vertex operator algebra $\mathcal{U}$.

Recall that we have seen $\alpha_{-}^{2}=\frac{p}{p'}$ in Section
\ref{subsec:The-unitary-series}. Now let $\alpha_{+}^{2}=\frac{p'}{p}$,
$x=\exp\left(2\pi i\alpha_{+}^{2}\right)$, $y=\exp\left(2\pi i\alpha_{-}^{2}\right)$,
$\left[l\right]=x^{l/2}-x^{-l/2}$, $\left[l'\right]=y{}^{l'/2}-y{}^{-l'/2}.$
Now we fix central charge $c_{p}$, denote $L\left(c_{p},h_{\left(i',i\right)}^{\left(p\right)}\right)$
by $\left(i',i\right)$. Let $\left(a',a\right)$, $\left(m',m\right)$,
$\left(n',n\right)$, $\left(c',c\right)$, $\left(b',b\right)$,
$\left(d',d\right)$ be irreducible $L\left(c_{p},0\right)$-modules,
the braiding matrices of screened vertex operators have the almost
factorized form (cf. (2.19) of \cite{FFK}):
\begin{alignat}{1}
 & \left(\tilde{B}_{\left(m',m\right),\left(n',n\right)}^{\left(a',a\right),\left(c',c\right)}\right){}_{\left(b',b\right),\left(d',d\right)}\nonumber \\
 & =i^{-\left(m'-1\right)\left(n-1\right)-\left(n'-1\right)\left(m-1\right)}\left(-1\right)^{1/2\left(a-b+c-d\right)\left(n'+m\right)+1/2\left(a'-b'+c'-d'\right)\left(n+m\right)}\label{FFK 2.19}\\
 & \cdot r\left(a',m',n',c'\right)_{b',d'}\cdot r\left(a,m,n,c\right)_{b,d},\nonumber
\end{alignat}
where the nonvanishing matrix elements of $r$-matrices are

\begin{gather}
r\left(a,1,n,c\right)_{a,c}=r\left(a,m,1,c\right)_{c,a}=1,\nonumber \\
r\left(l\pm2,2,2,l\right)_{l\pm1,l\pm1}=x^{1/4},\nonumber \\
r\left(l,2,2,l\right)_{l\pm1,l\pm1}=\mp x^{-1/4\mp l/2}\frac{\left[1\right]}{\left[l\right]},\nonumber \\
r\left(l,2,2,l\right)_{l\pm1,l\mp1}=x^{-1/4}\frac{\left[l\pm1\right]}{\left[l\right]},\label{FFK 2.20}
\end{gather}
and the other $r$-matrices are given by the recursive relation
\begin{gather}
r\left(a,m+1,n,c\right)_{b,d}=\sum_{d_{1}\ge1}r\left(a,2,n,d_{1}\right)_{a_{1},d}\cdot r\left(a_{1},m,n,c\right)_{b,d_{1}},\nonumber \\
r\left(a,m,n+1,c\right)_{b,d}=\sum_{d_{1}\ge1}r\left(a,m,2,c_{1}\right)_{b,d_{1}}\cdot r\left(d_{1},m,n,c\right)_{c_{1},d},\label{FFK 2.21}
\end{gather}
for any choice of $a_{1}$ and $c_{1}$ compatible with the fusion
rules. The $r'$ matrices are given by the same formulas with the
replacement $x\to x'$, $\left[\ \ \ \ \right]\to\left[\ \ \ \ \right]'.$

Now we consider braiding matrix for $L\left(\frac{25}{28},0\right)$-modules.
Denote irreducible $L\left(\frac{25}{28},0\right)$-modules $L\left(\frac{25}{28},\frac{34}{7}\right)$
and $L\left(\frac{25}{28},\frac{9}{7}\right)$ by $Q_{2}$ and $Q_{3}$
respectively. For convenience, we will denote $\left(\tilde{B}_{Q_{a},Q_{b}}^{Q_{c},Q_{d}}\right)_{Q_{e},Q_{f}}$
by $\left(\tilde{B}_{a,b}^{c,d}\right)_{e,f}$, $a,b,c$, $d,$ $e$,
$f$ $\in\left\{ 2,3\right\} $. Now we are ready to give the following
lemma.

\begin{lemma} \label{(2332)_3,2 nonzero} \label{(3322)_3,3 nonzero}
\label{(2333)_3,2 nonzero} \label{(3,2,3,3)_2,3 nonzero}$\left(\tilde{B}_{2,2}^{3,3}\right)_{3,2}\not=0$,
$\left(\tilde{B}_{3,2}^{3,2}\right)_{3,3}\not=0$, $\left(\tilde{B}_{2,3}^{3,3}\right)_{3,2}\not=0$,
and $\left(\tilde{B}_{3,3}^{2,3}\right)_{2,3}\not=0$. \end{lemma}

\begin{proof} Using (\ref{FFK 2.19}), (\ref{FFK 2.20}), and Remark
\ref{pairs corr. to U's weights for unitary model}, to prove $\left(\tilde{B}_{2,2}^{3,3}\right)_{3,2}\not=0$,
it suffices to show that $r\left(5,3,3,5\right)_{3,5}\not=0$. Using
(\ref{FFK 2.20}) and (\ref{FFK 2.21}) we obtain:

\[
r\left(5,3,3,5\right)_{3,5}=r\left(5,2,3,4\right)_{4,5}\cdot r\left(4,2,3,5\right)_{3,4}+r\left(5,2,3,6\right)_{4,5}\cdot r\left(4,2,3,5\right)_{3,6}
\]

with
\begin{alignat*}{1}
r\left(5,2,3,4\right)_{4,5} & =r\left(5,2,2,5\right)_{4,4}\cdot r\left(4,2,2,4\right)_{5,5}+r\left(5,2,2,5\right)_{4,6}\cdot r\left(6,2,2,4\right)_{5,5}\\
 & =\frac{\left[4\right]\left[4\right]-\left[1\right]\left[1\right]}{\left[4\right]\left[5\right]},
\end{alignat*}
\begin{alignat*}{1}
r\left(4,2,3,5\right)_{3,4} & =r\left(4,2,2,4\right)_{3,3}\cdot r\left(3,2,2,5\right)_{4,4}+r\left(4,2,2,4\right)_{3,5}\cdot r\left(5,2,2,5\right)_{4,4}\\
 & =x^{2}\left(\frac{\left[1\right]}{\left[4\right]}+\frac{\left[3\right]\left[1\right]}{\left[4\right]\left[5\right]}\right),
\end{alignat*}
\begin{alignat*}{1}
r\left(5,2,3,6\right)_{4,5} & =r\left(5,2,2,5\right)_{4,4}\cdot r\left(4,2,2,6\right)_{5,5}+r\left(5,2,2,5\right)_{4,6}\cdot r\left(6,2,2,6\right)_{5,5}\\
 & =x^{5/2}\frac{\left[1\right]\left[6\right]+\left[4\right]\left[1\right]}{\left[5\right]\left[6\right]},
\end{alignat*}

\[
r\left(4,2,3,5\right)_{3,6}=r\left(4,2,2,4\right)_{3,5}\cdot r\left(5,2,2,5\right)_{4,6}=x^{-1/2}\frac{\left[3\right]}{\left[5\right]},
\]
where $\left[l\right]=2i\sin\left(\frac{8}{7}\pi l\right)$, $x=\exp\left(\frac{16}{7}\pi i\right)$.
Direct computation gives:
\begin{align*}
 & r\left(5,3,3,5\right)_{3,5}\\
 & =x^{2}\cdot\left(\frac{\left[4\right]^{2}-\left[1\right]^{2}}{\left[4\right]\left[5\right]}\left(\frac{\left[1\right]}{\left[4\right]}+\frac{\left[3\right]\left[1\right]}{\left[4\right]\left[5\right]}\right)+\frac{\left[1\right]\left[6\right]+\left[4\right]\left[1\right]}{\left[5\right]\left[6\right]}\cdot\frac{\left[3\right]}{\left[5\right]}\right)\\
 & =x^{2}\cdot\left(1+2\sin\left(\frac{\pi}{14}\right)+2\cos\left(\frac{\pi}{7}\right)\right)\not=0
\end{align*}
and hence $\left(\tilde{B}_{2,2}^{3,3}\right)_{3,2}\not=0$.

Similarly, to prove $\left(\tilde{B}_{3,2}^{3,2}\right)_{3,3}\not=0$,
$\left(\tilde{B}_{2,3}^{3,3}\right)_{3,2}\not=0$, and $\left(\tilde{B}_{3,3}^{2,3}\right)_{2,3}\not=0$,
it suffices to show that $r\left(3,3,5,5\right)_{3,3}\not=0$, $r\left(5,3,3,3\right)_{3,5}\not=0,$
and $r\left(3,5,3,3\right)_{5,3}\not=0$ respectively. Direct calculation
gives:

\[
r\left(3,3,5,5\right)_{3,3}=\frac{x^{2}}{8}\cdot\sin\left(\frac{\pi}{7}\right)\sec^{2}\left(\frac{\pi}{14}\right)\sec^{3}\left(\frac{3\pi}{14}\right)\left(-1+\sin\frac{\pi}{14}-2\cos\left(\frac{\pi}{7}\right)\right)\not=0,
\]

\[
r\left(5,3,3,3\right)_{3,5}=-\sin\left(\frac{\pi}{7}\right)\sec^{3}\left(\frac{\pi}{14}\right)\sec\left(\frac{3\pi}{14}\right)\left(\sin\frac{\pi}{7}+\cos\frac{\pi}{14}\right)\left(\cos^{2}\left(\frac{\pi}{14}\right)-\sin^{2}\left(\frac{\pi}{7}\right)\right)\not=0,
\]
and
\begin{alignat*}{1}
 & r\left(3,5,3,3\right)_{5,3}\\
 & =x^{-2}\cdot\sin^{3}\left(\frac{\pi}{7}\right)\left(\cos\left(\frac{\pi}{14}\right)+\cos\left(\frac{3\pi}{14}\right)\right)\sec^{5}\left(\frac{\pi}{14}\right)\sec\left(\frac{3\pi}{14}\right)\\
 & \cdot\left(\cos^{2}\left(\frac{\pi}{14}\right)-\sin^{2}\left(\frac{\pi}{7}\right)\right)\left(\cos^{2}\left(\frac{3\pi}{14}\right)-\sin^{2}\left(\frac{\pi}{7}\right)\right)\\
 & +x^{-2}\cdot\left(-6\sin\left(\frac{\pi}{7}\right)+9\cos\left(\frac{\pi}{14}\right)-7\cos\left(\frac{3\pi}{14}\right)\right)/\left(7\sin\left(\frac{\pi}{7}\right)+12\cos\left(\frac{\pi}{14}\right)+11\cos\left(\frac{3\pi}{14}\right)\right)\\
 & \not=0.
\end{alignat*}
Therefore $\left(\tilde{B}_{3,2}^{3,2}\right)_{3,3}\not=0$, $\left(\tilde{B}_{3,3}^{2,3}\right)_{3,2}\not=0$,
and $\left(\tilde{B}_{3,3}^{2,3}\right)_{2,3}\not=0$. \end{proof}

\subsection{GKO construction of the unitary Virasoro VOA}

Let $e,f$ and $h$ be the generators of $\mathfrak{sl}_{2}(\mathbb{C})$
such that
\[
\left[e,f\right]=h,\ \left[h,e\right]=2e,\ \left[h,f\right]=-2f.
\]
Let $\left\langle \cdot,\cdot\right\rangle $ be the standard invariant
bilinear form on $\mathfrak{sl}_{2}(\mathbb{C})$ defined by
\[
\left\langle h,h\right\rangle =2,\ \left\langle e,f\right\rangle =1,\ \left\langle e,e\right\rangle =\left\langle f,f\right\rangle =\left\langle h,e\right\rangle =\left\langle h,f\right\rangle =0.
\]
Let $\hat{\mathfrak{sl}}_{2}\left(\mathbb{C}\right)$ be the corresponding
affine algebra of type $A_{1}^{(1)}$ and $\lambda_{0},\ \lambda_{1}$
the fundamental weights for $\hat{\mathfrak{sl}}_{2}\left(\mathbb{C}\right)$.
Denote
\[
\mathcal{L}(m,k)=\mathcal{L}\left(\left(m-k\right)\lambda_{0}+k\lambda_{1}\right)
\]
the irreducible highest weight module of $\hat{\mathfrak{sl}}_{2}\left(\mathbb{C}\right)$-module
with highest weight $\left(m-k\right)\lambda_{0}+k\lambda_{1}$. It
was proved in \cite{FZ} that $\mathcal{L}\left(m,0\right)$ has a
natural vertex operator algebra structure for $m\in\mathbb{Z}_{+}$.
The Virasoro vector $\omega^{m}$ of $\mathcal{L}\left(m,0\right)$
is given by
\[
\omega^{m}=\frac{1}{2\left(m+2\right)}\left(\frac{1}{2}h_{-1}h+e_{-1}f+f_{-1}e\right)
\]
with central charge $\frac{3m}{m+2}$. Let $m\in\mathbb{Z}_{+}$,
then $\mathcal{L}\left(m,0\right)$ is a rational vertex operator
algebra and $\left\{ \mathcal{L}\left(m,\ k\right)|k=0,1,\cdots,m\right\} $
is the set of all the irreducible $\mathcal{L}\left(m,0\right)$-modules.
Moreover, the fusion rules are given by

\[
\mathcal{L}\left(m,j\right)\boxtimes\mathcal{L}\left(m,k\right)=\sum_{i=\max\left\{ 0,j+k-m\right\} }^{\min\left\{ j,k\right\} }\mathcal{L}\left(m,j+k-2i\right).
\]

Let $\mathcal{L}\left(m,0\right)_{1}$ be the weight 1 subspace of
$\mathcal{L}(m,0)$. Then $\mathcal{L}(m,0)_{1}$ has a structure
of Lie algebra isomorphic to $\mathfrak{sl}_{2}(\mathbb{C})$ under
$[a,b]=a_{0}b$, $\forall a,b\in\mathcal{L}(m,0)$. Let $h^{m}$,
$e^{m}$, $f^{m}$ be the generators of $\mathfrak{sl}_{2}(\mathbb{C})$
in $\mathcal{L}(m,0)_{1}$. Then $h^{m+1}:=h^{1}\otimes1+1\otimes h^{m},$
$e^{m+1}:=e^{1}\otimes1+1\otimes e^{m}$ and $f^{m+1}:=f^{1}\otimes1+1\otimes f^{m}$
generate a vertex operator subalgebra isomorphic to $\mathcal{L}\left(m+1,0\right)$
in $\mathcal{L}\left(m,1\right)\otimes\mathcal{L}\left(m,0\right)$.
It was proved in \cite{DL} and \cite{KR} that $\Omega^{m}:=\omega^{1}\otimes1+1\otimes\omega^{m}-\omega^{m+1}$
also gives a Virasoro vector with central charge $c_{m+2}=1-6/\left(m+2\right)\left(m+3\right)$.
Furthermore, $\omega^{m+1}$ and $\Omega^{m}$ are mutually commutative
and $\Omega^{m}$ generates a simple Virasoro vertex operator algebra
$L\left(c_{m+2},0\right)$. Therefore $\mathcal{L}\left(1,0\right)\otimes\mathcal{L}\left(m,0\right)$
contains a vertex operator subalgebra isomorphic to $L\left(c_{m+2},0\right)\otimes\mathcal{L}\left(m+1,0\right)$.
Note that both $L\left(c_{m+2},0\right)$ and $\mathcal{L}\left(m+1,0\right)$
are rational and every $\mathcal{L}\left(1,0\right)\otimes\mathcal{L}\left(m,0\right)$-module
can be decomposed into irreducible $L(c_{m+2},0)\otimes\mathcal{L}(m+1,0)$-submodules.
We have the following decomposition \cite{GKO}:

\begin{equation}
\mathcal{L}\left(1,\epsilon\right)\otimes\mathcal{L}\left(m,n\right)=\bigoplus_{0\le s\le m+3,s\equiv n+\epsilon\ \mbox{mod}2}L\left(c_{m+2},h_{\left(s+1,n+1\right)}^{(m+2)}\right)\otimes\mathcal{L}\left(m+1,s\right)\label{eq:(2.3)}
\end{equation}
where $\epsilon=0,1$ and $0\le n\le m$. This is the GKO-construction
of the unitary Virasoro vertex operator algebras.

\section{Structure of the $6A$-algebra $\mathcal{U}$}

Certain coset subalgebra of $V_{\sqrt{2}E_{8}}$ associated with extended
$E_{8}$ diagram is constructed in \cite{LYY2} by removing one node
from the diagram. In each case, the coset subalgebra contains some
Ising vectors and the coset subalgebra is generated by two Ising vectors
with inner product the same as the number given in the table in Section
1. In particular, the coset subalgebra $\mathcal{U}$ corresponding
to the $6A$ case was constructed, i.e., the case with inner product
$\frac{5}{2^{10}}.$ Let $\mathcal{V}$ be the $3A$-algebra, that
is, the vertex operator algebra generated by two Ising vectors whose
$\tau$-involutions generate $S_{3}$ and with inner product $\frac{13}{2^{10}}$.
The candidates for $\mathcal{V}$ were given \cite{M2} and it was
proved in \cite{SY} that only one of these candidates actually exists
and that there is unique vertex operator algebra structure on it.
Actually

\begin{eqnarray*}
\mathcal{V} & \cong & \left(\left(L\left(\frac{4}{5},0\right)\oplus L\left(\frac{4}{5},3\right)\right)\otimes\left(L\left(\frac{6}{7},0\right)\oplus L\left(\frac{6}{7},5\right)\right)\right)\\
 &  & \oplus L\left(\frac{4}{5},\frac{2}{3}\right)^{+}\otimes L\left(\frac{6}{7},\frac{4}{3}\right)^{+}\oplus L\left(\frac{4}{5},\frac{2}{3}\right)^{-}\otimes L\left(\frac{6}{7},\frac{4}{3}\right)^{-}.
\end{eqnarray*}

Now we recall the following results about the $3A$-algebra $\mathcal{V}$
from \cite{SY} .

\begin{lemma} \label{rationality of 3A}The $3A$-algebra $\mathcal{V}$
is rational. \end{lemma}

\begin{lemma} \label{U_3A modules} All the irreducible $\mathcal{V}$-modules
are as follows:

\begin{align*}
\mathcal{V}=\mathcal{V}\left(0\right) & =\left(\left(L\left(\frac{4}{5},0\right)\oplus L\left(\frac{4}{5},3\right)\right)\otimes\left(L\left(\frac{6}{7},0\right)\oplus L\left(\frac{6}{7},5\right)\right)\right)\\
 & \oplus L\left(\frac{4}{5},\frac{2}{3}\right)^{+}\otimes L\left(\frac{6}{7},\frac{4}{3}\right)^{+}\oplus L\left(\frac{4}{5},\frac{2}{3}\right)^{-}\otimes L\left(\frac{6}{7},\frac{4}{3}\right)^{-},
\end{align*}

\begin{align*}
\mathcal{V}\left(\frac{1}{7}\right) & =\left(L\left(\frac{4}{5},0\right)\oplus L\left(\frac{4}{5},3\right)\right)\otimes\left((L\left(\frac{6}{7},\frac{1}{7}\right)\oplus L\left(\frac{6}{7},\frac{22}{7}\right)\right)\\
 & \oplus L(\left(\frac{4}{5},\frac{2}{3}\right)^{+}\otimes L\left(\frac{6}{7},\frac{10}{21}\right)^{+}\oplus L\left(\frac{4}{5},\frac{2}{3}\right)^{-}\otimes L\left(\frac{6}{7},\frac{10}{21}\right)^{-},
\end{align*}

\begin{align*}
\mathcal{V}\left(\frac{5}{7}\right) & =\left(L\left(\frac{4}{5},0\right)\oplus L\left(\frac{4}{5},3\right)\right)\otimes\left(L\left(\frac{6}{7},\frac{5}{7}\right)\oplus L\left(\frac{6}{7},\frac{12}{7}\right)\right)\\
 & \oplus L(\left(\frac{4}{5},\frac{2}{3}\right)^{+}\otimes L\left(\frac{6}{7},\frac{1}{21}\right)^{+}\oplus L\left(\frac{4}{5},\frac{2}{3}\right)^{-}\otimes L\left(\frac{6}{7},\frac{1}{21}\right)^{-},
\end{align*}

\begin{align*}
\mathcal{V}\left(\frac{2}{5}\right) & =\left(L\left(\frac{4}{5},\frac{2}{5}\right)\oplus L\left(\frac{4}{5},\frac{7}{5}\right)\right)\otimes\left(L\left(\frac{6}{7},0\right)\oplus L\left(\frac{6}{7},5\right)\right)\\
 & \oplus L\left(\frac{4}{5},\frac{1}{15}\right)^{+}\otimes L(\left(\frac{6}{7},\frac{4}{3}\right)^{+}\oplus L\left(\frac{4}{5},\frac{1}{15}\right)^{-}\otimes L\left(\frac{6}{7},\frac{4}{3}\right)^{-},
\end{align*}

\begin{align*}
\mathcal{V}\left(\frac{19}{35}\right) & =\left(L\left(\frac{4}{5},\frac{2}{5}\right)\oplus L\left(\frac{4}{5},\frac{7}{5}\right)\right)\otimes\left(L\left(\frac{6}{7},\frac{1}{7}\right)\oplus L\left(\frac{6}{7},\frac{22}{7}\right)\right)\\
 & \oplus L(\left(\frac{4}{5},\frac{1}{15}\right)^{+}\otimes L(\left(\frac{6}{7},\frac{10}{21}\right)^{+}\oplus L(\left(\frac{4}{5},\frac{1}{15}\right)^{-}\otimes L(\left(\frac{6}{7},\frac{10}{21}\right)^{-},
\end{align*}
\begin{align*}
\mathcal{V}\left(\frac{39}{35}\right) & =\left(L\left(\frac{4}{5},\frac{2}{5}\right)\oplus L\left(\frac{4}{5},\frac{7}{5}\right)\right)\otimes\left(L\left(\frac{6}{7},\frac{5}{7}\right)\oplus L\left(\frac{6}{7},\frac{12}{7}\right)\right)\\
 & \oplus L(\left(\frac{4}{5},\frac{1}{15}\right)^{+}\otimes L(\left(\frac{6}{7},\frac{1}{21}\right)^{+}\oplus L(\left(\frac{4}{5},\frac{1}{15}\right)^{-}\otimes L\left(\frac{6}{7},\frac{1}{21}\right)^{-}.
\end{align*}

\end{lemma}

\begin{proposition}\label{fusion rules for U3A} Fusion rules for
all the irreducible $\mathcal{V}$-modules are as the following. (For
simplicity, we denote $\mathcal{V}\left(k\right)$ by $k$, where
$k=0,\frac{1}{7},\frac{5}{7},\frac{2}{5},\frac{19}{35},\frac{39}{35}$.
)
\begin{center}
\begin{tabular}{|c|c|c|c|c|c|}
\hline
0  & $\frac{1}{7}$  & $\frac{5}{7}$  & $\frac{2}{5}$  & $\frac{19}{35}$  & $\frac{39}{35}$\tabularnewline
\hline
\hline
$\frac{1}{7}$  & $0+\frac{5}{7}$  & $\frac{1}{7}+\frac{5}{7}$  & $\frac{19}{35}$  & $\frac{2}{5}+\frac{39}{35}$  & $\frac{19}{35}+\frac{39}{35}$\tabularnewline
\hline
$\frac{5}{7}$  & $\frac{1}{7}+\frac{5}{7}$  & $0+\frac{1}{7}+\frac{5}{7}$  & $\frac{39}{35}$  & $\frac{19}{35}+\frac{39}{35}$  & $\frac{2}{5}+\frac{19}{35}+\frac{39}{35}$\tabularnewline
\hline
$\frac{2}{5}$  & $\frac{19}{35}$  & $\frac{39}{35}$  & $0+\frac{2}{5}$  & $\frac{1}{7}+\frac{39}{35}$  & $\frac{5}{7}+\frac{39}{35}$\tabularnewline
\hline
$\frac{19}{35}$  & $\frac{2}{5}+\frac{39}{35}$  & $\frac{19}{35}+\frac{39}{35}$  & $\frac{1}{7}+\frac{39}{35}$  & $0+\frac{5}{7}+\frac{2}{5}+\frac{39}{35}$  & $\frac{1}{7}+\frac{5}{7}+\frac{19}{35}+\frac{39}{35}$\tabularnewline
\hline
$\frac{39}{35}$  & $\frac{19}{35}+\frac{39}{35}$  & $\frac{2}{5}+\frac{19}{35}+\frac{39}{35}$  & $\frac{5}{7}+\frac{39}{35}$  & $\frac{1}{7}+\frac{5}{7}+\frac{19}{35}+\frac{39}{35}$  & $0+\frac{1}{7}+\frac{5}{7}+\frac{2}{5}+\frac{19}{35}+\frac{39}{35}$\tabularnewline
\hline
\end{tabular}
\par\end{center}

\end{proposition}

It was proved in \cite{LYY2} that $\mathcal{V}\subset\mathcal{U}$
and as a module of $\mathcal{V}\otimes L\left(\frac{25}{28},0\right)$,

\[
\mathcal{U}\cong\mathcal{V}\otimes L\left(\frac{25}{28},0\right)\oplus\mathcal{V}\left(\frac{1}{7}\right)\otimes L\left(\frac{25}{28},\frac{34}{7}\right)\oplus\mathcal{V}\left(\frac{5}{7}\right)\otimes L\left(\frac{25}{28},\frac{9}{7}\right).
\]
From here forward, we denote
\begin{gather}
P_{1}=\mathcal{V},\ P_{2}=\mathcal{V}\left(\frac{1}{7}\right),\ P_{3}=\mathcal{V}\left(\frac{5}{7}\right),\nonumber \\
Q_{1}=L\left(\frac{25}{28},0\right),\ Q_{2}=L\left(\frac{25}{28},\frac{34}{7}\right),\ Q_{3}=L\left(\frac{25}{28},\frac{9}{7}\right),\label{PQ notation}
\end{gather}
and $U^{i}=P_{i}\otimes Q_{i}$, $i=1,2,3$. Then
\[
\mathcal{U}\cong P_{1}\otimes Q_{1}\oplus P_{2}\otimes Q_{2}\oplus P_{3}\otimes Q_{3}=U^{1}\oplus U^{2}\oplus U^{3}.
\]

\begin{remark} \label{rationality-1} Since $\mathcal{V}\otimes L\left(\frac{25}{28},0\right)$
is a rational and $C_{2}$-cofinite vertex operator algebra, it is
straightforward to see that $\mathcal{U}$ is also rational and $C_{2}$-cofinite
by \cite{HKL,ABD}. \end{remark}

\begin{remark} \label{self-dual} 1. By fusion rules for irreducible
$L\left(\frac{25}{28},0\right)$-modules and $\mathcal{V}$-modules
in Propositions \ref{fusion rules of virasoro modules} and \ref{fusion rules for U3A},
and rationality of the $3A$-algebra in Lemma \ref{rationality of 3A},
we see that to study $\mathcal{U}$, we shall study an extension of
a rational vertex operator algebra $U^{1}$ by two $U^{1}$-irreducible
modules $U^{2},U^{3}$ which are not simple current modules.

2. Since $\mathcal{U}_{1}=0$ and $\dim\mathcal{U}_{0}=1$ by Theorem
\ref{bilinear form}, there is a unique bilinear form on $\mathcal{U}$
and thus $\mathcal{U}'\cong\mathcal{U}$. Without loss of generality,
we can identify $\mathcal{U}$ with $\mathcal{U}'$.\end{remark}

\begin{remark}Let $h_{1},h_{2},h_{3},h_{4}$ be four irreducible
$\mathcal{V}$-modules, and fix a basis of intertwining operators.
By Section \ref{subsec:Braiding-matrix},  there exists $\left(B_{_{h_{4},h_{1}}}^{h_{3,}h_{2}}\right)_{\mu,\gamma}^{i,j;k,l}\in\mathbb{C}$
such that
\begin{alignat}{1}
 & E\left\langle u_{h_{4}'},\mathcal{Y}_{h_{3},\mu;i}^{h_{4}}\left(u_{h_{3}},z_{1}\right)\mathcal{Y}_{h_{2},h_{1};j}^{\mu}\left(u_{h_{2}},z_{2}\right)u_{h_{1}}\right\rangle \nonumber \\
 & =\sum_{k,l,\gamma}\left(B_{_{h_{4},h_{1}}}^{h_{3,}h_{2}}\right)_{\mu,\gamma}^{i,j;k,l}E\left\langle u_{h_{4}'},\mathcal{Y}_{h_{2},\gamma;k}^{h_{4}}\left(u_{h_{2}},z_{1}\right)\mathcal{Y}_{h_{3},h_{1};l}^{\gamma}\left(u_{h_{3}},z_{2}\right)u_{h_{1}}\right\rangle \label{B}
\end{alignat}

\end{remark}

\subsection{Uniqueness of VOA structure on $\mathcal{U}$}

Recall notations in (\ref{PQ notation}). For convenience, we list
fusion rules of $I_{P_{1}}\left(_{P_{a}\ P_{b}}^{P_{c}}\right)$ and
$I_{Q_{1}}\left(_{Q_{a}\ Q_{b}}^{Q_{c}}\right)$ with $a,b,c\in\left\{ 1,2,3\right\} $
from Propositions \ref{fusion rules of virasoro modules} and \ref{fusion rules for U3A}
in the following table.
\begin{center}
\begin{tabular}{|c|c|c|}
\hline
$P_{1}$  & $P_{2}$  & $P_{3}$\tabularnewline
\hline
\hline
$P_{2}$  & $P_{1}+P_{3}$  & $P_{2}$+$P_{3}$\tabularnewline
\hline
$P_{3}$  & $P_{2}+P_{3}$  & $P_{1}+P_{2}$+$P_{3}$\tabularnewline
\hline
\end{tabular}\label{Fusion rules for P_i's}
\par\end{center}

\begin{center}
\begin{tabular}{|c|c|c|}
\hline
$Q_{1}$  & $Q_{2}$  & $Q_{3}$\tabularnewline
\hline
\hline
$Q_{2}$  & $Q_{1}+Q_{3}$  & $Q_{2}+Q_{3}$\tabularnewline
\hline
$Q_{3}$  & $Q_{2}+Q_{3}$  & $Q_{1}+Q_{2}+Q_{3}$\tabularnewline
\hline
\end{tabular}\label{Fusion rules for Q_i's}
\par\end{center}

Note that the fusion rules $N_{P_{1}}\left(_{P_{a}\ P_{b}}^{P_{c}}\right)=N_{Q_{1}}\left(_{Q_{a}\ Q_{b}}^{Q_{c}}\right)$,
which is either $0$ or $1$. We immediately get
\[
N_{a,b}^{c}=N_{U^{1}}\left(_{U^{a}\ U^{b}}^{U^{c}}\right)=N_{P_{1}}\left(_{P_{a}\ P_{b}}^{P_{c}}\right)\cdot N_{Q_{1}}\left(_{Q_{a}\ Q_{b}}^{Q_{c}}\right)=N_{P_{1}}\left(_{P_{a}\ P_{b}}^{P_{c}}\right)=N_{Q_{1}}\left(_{Q_{a}\ Q_{b}}^{Q_{c}}\right).
\]

We fix a basis $\overline{\mathcal{Y}}_{a,b}^{c}\in I_{Q_{1}}\left(_{Q_{a}\ Q_{b}}^{Q_{c}}\right)$
as in Section \ref{subsec:Braiding-matrix} and \cite{FFK}, and choose
an arbitrary basis of $\mathcal{Y}_{a,b}^{c}\in I_{P_{1}}\left(_{P_{a}\ P_{b}}^{P_{c}}\right).$
Then $\mathcal{I}_{a,b}^{c}=\mathcal{Y}_{a,b}^{c}\otimes\overline{\mathcal{Y}}_{a,b}^{c}$
is a basis of $I_{U^{1}}\left(_{U^{a}\ U^{b}}^{U^{c}}\right)$.

Now let $\left(\mathcal{U},Y\right)$ be a vertex operator algebra
structure on $\mathcal{U}$ with

\[
Y\left(u_{1}^{a}\otimes u_{2}^{a},z\right)=\sum_{b,c\in\left\{ 1,2,3\right\} }\text{\ensuremath{\lambda}}_{a,b}^{c}\cdot\mathcal{I}_{a,b}^{c}\left(u_{1}^{a}\otimes u_{2}^{a},z\right)=\sum_{b,c\in\left\{ 1,2,3\right\} }\text{\ensuremath{\lambda}}_{a,b}^{c}\mathcal{Y}_{a,b}^{c}(u_{1}^{a},z)\otimes\overline{\mathcal{Y}}_{a,b}^{c}(u_{2}^{b},z)
\]
where $u_{1}^{a}\in P^{a},u_{2}^{a}\in Q^{a}$.

The following lemma plays an important role in the proof of the uniqueness
of the vertex operator algebra structure on $\mathcal{U}$.

\begin{lemma} $\text{\ensuremath{\lambda}}_{a,b}^{c}\neq0$ if $I_{U^{1}}\left(_{U^{a}\ U^{b}}^{U^{c}}\right)\neq0$.
\end{lemma}
\begin{claim}
$\lambda_{k,1}^{k}\not=0,\forall k=2,3$.
\end{claim}
\begin{proof}
For any $u^{k}\in U^{k}$, $k=1,2,3$, using skew symmetry of $Y\left(\cdot,z\right)$
(\cite{FHL}), we have

\[
Y(u^{k},z)u^{1}=e^{zL\left(-1\right)}Y\left(u^{1},-z\right)u^{k}=\lambda_{1,k}^{k}\cdot e^{zL\left(-1\right)}\mathcal{I}_{1,k}^{k}\left(u^{1},-z\right)u^{k}=\lambda_{k,1}^{k}\cdot\mathcal{I}_{k,1}^{k}\left(u^{k},z\right)u^{1}.
\]
Since $U^{k}$ is an irreducible $U^{1}$-module, we have $\lambda_{1,k}^{k}\not=0$,
$\forall k=1,2,3$. So $\lambda_{k,1}^{k}\not=0$, $\forall k=1,2,3.$
\end{proof}
\begin{claim}
$\lambda_{k,k}^{1}\not=0$, $\forall k=2,3$.
\end{claim}
\begin{proof}
Note that from Remark \ref{self-dual}, $\mathcal{U}$ has a unique
invariant bilinear form $\left\langle \cdot,\cdot\right\rangle $
with $\left\langle 1,1\right\rangle =1$. For $u^{k},v^{k}\in U^{k},$
$k=1,2,3$, we have
\[
\left\langle Y\left(u^{k},z)v^{k}\right),u^{1}\right\rangle =\left\langle v^{k},Y\left(e^{zL\left(-1\right)}\left(-z^{-2}\right)^{L\left(0\right)}u^{k},z^{-1}\right)u^{1}\right\rangle .
\]
That is,
\[
\left\langle \lambda_{k,k}^{1}\cdot\mathcal{I}_{k,k}^{1}\left(u^{k},z\right)v^{k},u^{1}\right\rangle =\left\langle v^{k},\lambda_{k,1}^{k}\cdot\mathcal{I}_{k,1}^{k}\left(e^{zL\left(-1\right)}\left(-z^{-2}\right)^{L\left(0\right)}u^{k},z^{-1}\right)u^{1}\right\rangle .
\]
Applying previous claim, $\lambda_{k,1}^{k}\not=0$, and hence $\lambda_{k,k}^{1}\not=0$,
$\forall k=2,3.$
\end{proof}
\begin{claim}
\label{lambda_2,3,2; 2,3,3; 3,3,2; 2,2,3; nonzero}\emph{ $\lambda_{2,2}^{3}\not=0$,
$\lambda_{3,3}^{2}\not=0$, $\lambda_{2,3}^{k}\not=0,\lambda_{3,2}^{k}\not=0,\forall k=2,3$.
}
\end{claim}
\begin{proof}
Let $u^{2},v^{2}\in U^{2},$ $u^{3}\in U^{3}$. Skew symmetry of $Y$
gives
\[
\left\langle Y\left(u^{2},z\right)u^{3},v^{2}\right\rangle =\left\langle e^{zL\left(-1\right)}Y\left(u^{3},-z\right)u^{2},v^{2}\right\rangle ,
\]
that is,
\begin{equation}
\left\langle \lambda_{2,3}^{2}\cdot\mathcal{I}_{2,3}^{2}\left(u^{2},z\right)u^{3},v^{2}\right\rangle =\left\langle \lambda_{3,2}^{2}\cdot e^{zL\left(-1\right)}\mathcal{I}_{3,2}^{2}\left(u^{3},-z\right)u^{2},v^{2}\right\rangle .\label{nonzero-1}
\end{equation}
So $\lambda_{2,3}^{2}$ and $\lambda_{3,2}^{2}$ are both zero or
nonzero.

For any $u^{1}\in U^{1},$ $u^{2},v^{2},w^{2}\in U^{2}$ and $u^{3}\in U^{3}$,
commutativity of $Y$ in (\ref{communitivity}) implies

\[
\iota_{12}^{-1}\left\langle u^{1},Y\left(u^{2},z_{1}\right)Y\left(u^{3},z_{2}\right)v^{2}\right\rangle =\iota_{21}^{-1}\left\langle u^{1},Y\left(u^{3},z_{1}\right)Y\left(u^{2},z_{1}\right)v^{2}\right\rangle .
\]
That is,

\begin{alignat}{1}
\iota_{12}^{-1} & \left\langle u^{1},\lambda_{2,2}^{1}\lambda_{3,2}^{2}\cdot\mathcal{I}_{2,2}^{1}\left(u^{2},z_{1}\right)\mathcal{I}_{3,2}^{2}\left(u^{3},z_{2}\right)v^{2}\right\rangle \nonumber \\
 & =\iota_{21}^{-1}\left\langle u^{1},\lambda_{3,3}^{1}\lambda_{2,2}^{3}\cdot\mathcal{I}_{3,3}^{1}\left(u^{3},z_{2}\right)\mathcal{I}_{2,2}^{3}\left(u^{2},z_{1}\right)v^{2}\right\rangle ,\label{nonzero-2}
\end{alignat}
\begin{alignat}{1}
 & \iota_{12}^{-1}\left\langle w^{2},\lambda_{2,3}^{2}\lambda_{3,2}^{3}\cdot\mathcal{I}_{2,3}^{2}\left(u^{2},z_{1}\right)\mathcal{I}_{3,2}^{3}\left(u^{3},z_{2}\right)v^{2}\right\rangle \nonumber \\
 & =\iota_{21}^{-1}\left\langle w^{2},\lambda_{3,3}^{2}\lambda_{2,2}^{3}\cdot\mathcal{I}_{3,3}^{2}\left(u^{3},z_{2}\right)\mathcal{I}_{2,2}^{3}\left(u^{2},z_{1}\right)v^{2}\right\rangle .\label{nonzero-3}
\end{alignat}
Using (\ref{nonzero-2}), (\ref{nonzero-3}) and previous claim, we
see that either $\lambda_{2,3}^{2}=\lambda_{3,2}^{2}=\lambda_{2,2}^{3}=0$
or none of $\lambda_{2,3}^{2},$ $\lambda_{3,2}^{2},$ $\lambda_{2,2}^{3}$
is zero. For $i,j\in\left\{ 1,2,3\right\} $, denote $U^{i}.U^{j}=\left\langle u_{n}^{i}u^{j}|u^{i}\in U^{i},u^{j}\in U^{j},n\in\mathbb{Z}\right\rangle .$
Assume $\lambda_{2,3}^{2}=\lambda_{3,2}^{2}=\lambda_{2,2}^{3}=0,$
then we have $\left(U^{1}+U^{2}\right).\left(U^{1}+U^{2}\right)\subset U^{1}+U^{2}$,
$\left(U^{1}+U^{2}\right).U^{3}\subset U^{3}$, so $U^{1}+U^{2}$
is a vertex operator subalgebra of $\mathcal{U}$ and $U^{3}$ is
a $U^{1}+U^{2}$-module. Now $U^{1}+U^{2}$ is an extension of a ''good''
vertex operator algebra, so $U^{1}+U^{2}$ is rational by Theorem
\ref{rationality of extesnion of VOA }. Note that $U^{1}.U^{1}=U^{1}$,
$U^{1}.U^{2}=U^{2}$, $U^{2}.U^{1}=U^{2},$ $U^{2}.U^{2}=U^{1}$.
Define $\sigma:U^{1}+U^{2}\to U^{1}+U^{2}$ such that $\sigma|_{U^{1}}=1$
and $\sigma|_{U^{2}}=-1$. Then $\sigma$ is an order 2 automorphism
of $U^{1}+U^{2}$ with $\left(U^{1}+U^{2}\right)^{\sigma}=U^{1}$
and $U^{2}$ is a $U^{1}$-module. Apply quantum Galois theory in
Theorems \ref{classical galois theory} and \ref{quantum dimension and orbifold module},
$U^{2}$ is a simple current $U^{1}$-module, which is a contradiction.
Therefore, $\lambda_{2,3}^{2}\not=0$, $\lambda_{3,2}^{2}\not=0$
and $\lambda_{2,2}^{3}\not=0$.

Similarly, when $u^{k},v^{k},w^{k}\in U^{k}$, (\ref{communitivity})
gives
\begin{alignat*}{1}
 & \iota_{12}^{-1}\left\langle w^{2},\lambda_{2,3}^{2}\lambda_{3,2}^{3}\cdot\mathcal{I}_{2,3}^{2}\left(u^{2},z_{1}\right)\mathcal{I}_{3,2}^{3}\left(u^{3},z_{2}\right)v^{2}\right\rangle \\
 & =\iota_{21}^{-1}\left\langle w^{2},\lambda_{3,3}^{2}\lambda_{2,2}^{3}\cdot\mathcal{I}_{3,3}^{2}\left(u^{3},z_{2}\right)\mathcal{I}_{2,2}^{3}\left(u^{2},z_{1}\right)v^{2}\right\rangle .
\end{alignat*}
Since $\lambda_{2,3}^{2}\not=0$ and $\lambda_{2,2}^{3}\not=0$, we
see that either $\lambda_{3,2}^{3}=\lambda_{3,3}^{2}=0,$ or both
$\lambda_{3,2}^{3}$ and $\lambda_{3,3}^{2}$ are nonzero.

Assume $\lambda_{3,2}^{3}=\lambda_{3,3}^{2}=0,$ then by skew symmetry
of the vertex operator $Y$, we have $\lambda_{3,2}^{3}=\lambda_{3,3}^{2}=\lambda_{2,3}^{3}=0$.
Now $\left(U^{1}+U^{3}\right).\left(U^{1}+U^{3}\right)=U^{1}+U^{3}$
and $\left(U^{1}+U^{3}\right).U^{2}=U^{2},$ so $U^{1}+U^{3}$ is
a vertex operator algebra and $U^{2}$ is a $U^{1}+U^{3}$-module.
Also note that $U^{2}.U^{2}=U^{1}+U^{3},U^{2}.\left(U^{1}+U^{3}\right)=U^{2},$
so $U^{2}$ is a simple current module of $U^{1}+U^{3}$, which implies
$q\dim_{U^{1}}\left(U^{1}+U^{3}\right)=q\dim_{U^{1}}U^{2}$, i.e.,
\begin{equation}
1+q\dim_{U^{1}}U^{3}=q\dim_{U^{1}}U^{2}.\label{nonzero-prop.}
\end{equation}
Recall the fusion rules listed in Section \ref{Fusion rules for P_i's}
and the results of quantum dimensions (see Proposition \ref{product property of quantum dimension }).
For $i\in\{1,2,3\}$ we have
\begin{equation}
q\dim_{P_{1}}P_{i}=q\dim_{Q_{1}}Q_{i}\label{Pi-equals-Qi}
\end{equation}
\begin{equation}
q\dim_{P_{1}}P_{2}\cdot q\dim_{P_{1}}P_{2}=1+q\dim_{P_{1}}P_{3}\label{P2-fusion}
\end{equation}
Equation (\ref{nonzero-prop.}) and the equations above implies
\begin{equation}
1+q\dim_{P_{1}}P_{3}\cdot q\dim_{Q_{1}}Q_{3}=q\dim_{P_{1}}P_{2}\cdot q\dim_{Q_{1}}Q_{2}\label{Ui-fusion}
\end{equation}

Let $q\dim_{P_{1}}P_{2}=x>0$ and $q\dim_{P_{1}}P_{3}=y>0,$ we have
\[
\begin{cases}
1+y=x^{2}\\
1+y^{2}=x^{2}.
\end{cases}
\]
The previous system of equations holds if and only if $y=1.$ This
contradicts with that $P_{3}$ is not a simple current module of $P_{1}$.
Contradiction implies that $\lambda_{3,2}^{3}$, $\lambda_{2,3}^{3}$
, $\lambda_{3,3}^{2}$ are all nonzero.
\begin{claim}
$\lambda_{3,3}^{3}\not=0$.
\end{claim}
Fix a basis $\left\{ \overline{\mathcal{Y}}_{a,b;i}^{c}|i=1,\cdots,N_{Q_{a},Q_{b}}^{Q_{c}}\right\} $
for $I\left(_{Q_{a},Q_{b}}^{Q_{c}}\right)$, $a,b,c\in\left\{ 1,2,3\right\} $
as in \cite{FFK}. Consider the four point functions on $\left(U^{3},U^{2},U^{3},U^{3}\right).$
Let $B_{3,3}^{2,3}$ be as defined in (\ref{B}). Let $w_{1}^{2}\otimes w_{2}^{2}\in U^{2},$
$t_{1}^{3}\otimes t_{2}^{3},$ $p_{1}^{3}\otimes p_{2}^{3}$$,$ $u_{1}^{3}\otimes u_{2}^{3}\in U^{3}$,
we have
\begin{alignat}{1}
 & E\left\langle t_{1}^{3}\otimes t_{2}^{3},Y\left(w_{1}^{2}\otimes w_{2}^{2},z_{1}\right)Y\left(u_{1}^{3}\otimes u_{2}^{3},z_{2}\right)p_{1}^{3}\otimes p_{2}^{3}\right\rangle \nonumber \\
 & =E\langle t_{1}^{3}\otimes t_{2}^{3},\lambda_{2,2}^{3}\lambda_{3,3}^{2}\cdot\mathcal{Y}_{2,2}^{3}\otimes\overline{\mathcal{Y}}_{2,2}^{3}\left(w_{1}^{2}\otimes w_{2}^{2},z_{1}\right)\cdot\mathcal{Y}_{3,3}^{2}\otimes\overline{\mathcal{Y}}_{3,3}^{2}\left(u_{1}^{3}\otimes u_{2}^{3},z_{2}\right)\cdot p_{1}^{3}\otimes p_{2}^{3}\nonumber \\
 & +\lambda_{2,3}^{3}\lambda_{3,3}^{3}\cdot\mathcal{Y}_{2,2}^{3}\otimes\overline{\mathcal{Y}}_{2,2}^{3}\left(w_{1}^{2}\otimes w_{2}^{2},z_{1}\right)\cdot\mathcal{Y}_{3,3}^{3}\otimes\overline{\mathcal{Y}}_{3,3}^{3}\left(u_{1}^{3}\otimes u_{2}^{3},z_{2}\right)\cdot p_{1}^{3}\otimes p_{2}^{3}\rangle\nonumber \\
 & =E\langle t_{1}^{3}\otimes t_{2}^{3},\lambda_{2,2}^{3}\lambda_{3,3}^{2}\cdot\mathcal{Y}_{2,2}^{3}\left(w_{1}^{2},z_{1}\right)\mathcal{Y}_{3,3}^{2}\left(u_{1}^{3},z_{2}\right)p_{1}^{3}\otimes\overline{\mathcal{Y}}_{2,2}^{3}\left(w_{2}^{2},z_{1}\right)\overline{\mathcal{Y}}_{3,3}^{2}\left(u_{2}^{3},z_{2}\right)p_{2}^{3}\nonumber \\
 & +\lambda_{2,3}^{3}\lambda_{3,3}^{3}\cdot\mathcal{Y}_{2,3}^{3}\left(w_{1}^{2},z_{1}\right)\mathcal{Y}_{3,3}^{3}\left(u_{1}^{3},z_{2}\right)p_{1}^{3}\otimes\overline{\mathcal{Y}}_{2,3}^{3}\left(w_{2}^{2},z_{1}\right)\overline{\mathcal{Y}}_{3,3}^{3}\left(u_{2}^{3},z_{2}\right)p_{2}^{3}\nonumber \\
 & =E\langle t_{1}^{3}\otimes t_{2}^{3},\lambda_{2,2}^{3}\lambda_{3,3}^{2}\cdot\sum_{i=2,3}\left(B_{3,3}^{2,3}\right)_{2,i}\mathcal{Y}_{3,i}^{3}\left(u_{1}^{3},z_{2}\right)\mathcal{Y}_{2,3}^{i}\left(w_{1}^{2},z_{1}\right)p_{1}^{3}\nonumber \\
 & \otimes\sum_{j=2,3}\left(\tilde{B}_{3,3}^{2,3}\right)_{2,j}\mathcal{\overline{Y}}_{3,j}^{3}\left(u_{2}^{3},z_{2}\right)\mathcal{\overline{Y}}_{2,3}^{j}\left(w_{2}^{2},z_{1}\right)p_{2}^{3}\nonumber \\
 & +\lambda_{2,3}^{3}\lambda_{3,3}^{3}\cdot\sum_{i=2,3}\left(B_{3,3}^{2,3}\right)_{3,i}\mathcal{Y}_{3,i}^{3}\left(u_{1}^{3},z_{2}\right)\mathcal{Y}_{2,3}^{i}\left(w_{1}^{2},z_{1}\right)p_{1}^{3}\nonumber \\
 & \otimes\sum_{j=2,3}\left(\tilde{B}_{3,3}^{2,3}\right)_{3,j}\mathcal{Y}_{3,j}^{3}\left(u_{2}^{3},z_{2}\right)\mathcal{Y}_{2,3}^{j}\left(w_{2}^{2},z_{1}\right)p_{2}^{3}\rangle\label{(3,5,3,3)_2,3 nonzero-1}
\end{alignat}
In the mean time, we have
\begin{alignat}{1}
 & E\left\langle t_{1}^{3}\otimes t_{2}^{3},Y\left(u_{1}^{3}\otimes u_{2}^{3},z_{2}\right)Y\left(w_{1}^{2}\otimes w_{2}^{2},z_{1}\right)p_{1}^{3}\otimes p_{2}^{3}\right\rangle \nonumber \\
 & =E\langle t_{1}^{3}\otimes t_{2}^{3},\lambda_{3,2}^{3}\lambda_{2,3}^{2}\cdot\mathcal{Y}_{3,2}^{3}\left(u_{1}^{3},z_{2}\right)\mathcal{Y}_{2,3}^{2}\left(w_{1}^{2},z_{1}\right)p_{1}^{3}\otimes\mathcal{\overline{Y}}_{3,2}^{3}\left(u_{2}^{3},z_{2}\right)\mathcal{\overline{Y}}_{2,3}^{2}\left(w_{2}^{2},z_{1}\right)p_{2}^{3}\nonumber \\
 & +\lambda_{3,3}^{3}\lambda_{2,3}^{3}\cdot\mathcal{Y}_{3,3}^{3}\left(u_{1}^{3},z_{2}\right)\mathcal{Y}_{2,3}^{3}\left(w_{1}^{2},z_{1}\right)p_{1}^{3}\otimes\mathcal{\overline{Y}}_{3,3}^{3}\left(u_{2}^{3},z_{2}\right)\mathcal{\overline{Y}}_{2,3}^{3}\left(w_{2}^{2},z_{1}\right)p_{2}^{3}\label{(3,5,3,3)_2,3 nonzero-2}
\end{alignat}

(\ref{(3,5,3,3)_2,3 nonzero-1}) and (\ref{(3,5,3,3)_2,3 nonzero-2})
together with the linear independence of the four point functions
as mentioned in Section \ref{subsec:Braiding-matrix} imply that
\[
\begin{cases}
\lambda_{2,2}^{3}\lambda_{3,3}^{2}\left(B_{3,3}^{2,3}\right)_{2,2}\cdot\left(\tilde{B}_{3,3}^{2,3}\right)_{2,2}+\lambda_{2,3}^{3}\lambda_{3,3}^{3}\cdot\left(B_{3,3}^{2,3}\right)_{3,2}\cdot\left(\tilde{B}_{3,3}^{2,3}\right)_{3,2}=\lambda_{3,2}^{3}\lambda_{2,3}^{2}\\
\lambda_{2,2}^{3}\lambda_{3,3}^{2}\cdot\left(B_{3,3}^{2,3}\right)_{2,3}\cdot\left(\tilde{B}_{3,3}^{2,3}\right)_{2,3}+\lambda_{2,3}^{3}\lambda_{3,3}^{3}\cdot\left(B_{3,3}^{2,3}\right)_{3,3}\cdot\left(\tilde{B}_{3,3}^{2,3}\right)_{3,3}=\lambda_{3,3}^{3}\lambda_{2,3}^{3}\\
\lambda_{2,2}^{3}\lambda_{3,3}^{2}\left(B_{3,3}^{2,3}\right)_{2,2}\cdot\left(\tilde{B}_{3,3}^{2,3}\right)_{2,3}+\lambda_{2,3}^{3}\lambda_{3,3}^{3}\cdot\left(B_{3,3}^{2,3}\right)_{3,2}\cdot\left(\tilde{B}_{3,3}^{2,3}\right)_{3,3}=0\\
\lambda_{2,2}^{3}\lambda_{3,3}^{2}\cdot\left(B_{3,3}^{2,3}\right)_{2,3}\cdot\left(\tilde{B}_{3,3}^{2,3}\right)_{2,2}+\lambda_{2,3}^{3}\lambda_{3,3}^{3}\cdot\left(B_{3,3}^{2,3}\right)_{3,3}\cdot\left(\tilde{B}_{3,3}^{2,3}\right)_{3,2}=0.
\end{cases}
\]
Assume that $\lambda_{3,3}^{3}=0$. Then the above system of equations
become
\[
\begin{cases}
\lambda_{2,2}^{3}\lambda_{3,3}^{2}\left(B_{3,3}^{2,3}\right)_{2,2}\cdot\left(\tilde{B}_{3,3}^{2,3}\right)_{2,2}=\lambda_{3,2}^{3}\lambda_{2,3}^{2}\\
\lambda_{2,2}^{3}\lambda_{3,3}^{2}\cdot\left(B_{3,3}^{2,3}\right)_{2,3}\cdot\left(\tilde{B}_{3,3}^{2,3}\right)_{2,3}=0\\
\lambda_{2,2}^{3}\lambda_{3,3}^{2}\left(B_{3,3}^{2,3}\right)_{2,2}\cdot\left(\tilde{B}_{3,3}^{2,3}\right)_{2,3}=0\\
\lambda_{2,2}^{3}\lambda_{3,3}^{2}\cdot\left(B_{3,3}^{2,3}\right)_{2,3}\cdot\left(\tilde{B}_{3,3}^{2,3}\right)_{2,2}=0.
\end{cases}
\]
Since we already have proved $\lambda_{2,2}^{3}\not=0$, $\lambda_{3,2}^{3}\not=0$,
$\lambda_{2,3}^{2}\not=0$ and $\lambda_{3,3}^{2}\not=0$ in Claim
\ref{lambda_2,3,2; 2,3,3; 3,3,2; 2,2,3; nonzero}, the above system
of equations implies
\[
\begin{cases}
\left(B_{3,3}^{2,3}\right)_{2,2}\cdot\left(\tilde{B}_{3,3}^{2,3}\right)_{2,2}\not=0\\
\left(B_{3,3}^{2,3}\right)_{2,3}\cdot\left(\tilde{B}_{3,3}^{2,3}\right)_{2,3}=0\\
\left(B_{3,3}^{2,3}\right)_{2,2}\cdot\left(\tilde{B}_{3,3}^{2,3}\right)_{2,3}=0\\
\left(B_{3,3}^{2,3}\right)_{2,3}\cdot\left(\tilde{B}_{3,3}^{2,3}\right)_{2,2}=0.
\end{cases}
\]
By Lemma \ref{(3,2,3,3)_2,3 nonzero}, $\left(\tilde{B}_{3,3}^{2,3}\right)_{2,3}\not=0$.
So the third equation of the above system implies $\left(B_{3,3}^{2,3}\right)_{2,2}=0$,
which contradicts with the first equation of the above system. Contradiction
implies $\lambda_{3,3}^{3}\not=0$.
\end{proof}
Let $\left(\mathcal{U},\ Y\right)$ be a vertex operator algebra structure
on $\mathcal{U}$. First we fix a basis $\left\{ \overline{\mathcal{Y}}_{a,b}^{c}\left(\cdot,z\right)|a,b,c=1,2,3\right\} $
for space of intertwining operators of type $\left(_{Q_{a},Q_{b}}^{Q_{c}}\right)$,
$a,b,c\in\left\{ 1,2,3\right\} $ as in \cite{FFK}. Without loss
of generality, we can choose a basis $\left\{ \mathcal{Y}\left(\cdot,z\right)|a,b,c=1,2,3\right\} $
for space of intertwining operators of type $\left(_{P_{a},P_{b}}^{P_{c}}\right)$,
$a,b,c\in\left\{ 1,2,3\right\} $ such that the coefficients $\lambda_{a,b}^{c}=1$
if $N_{a,b}^{c}\not=0$. Fix $\mathcal{I}_{a,b}^{c}\left(\cdot,z\right)=\mathcal{Y}_{a,b}^{c}\left(\cdot,z\right)\otimes\overline{\mathcal{Y}}_{a,b}^{c}\left(\cdot,z\right).$Now
we have $\left(\mathcal{U},\ Y\right)$, a vertex operator algebra
structure on $\mathcal{U}=U^{1}\oplus U^{2}\oplus U^{3}$ such that
for any $u^{k},v^{k}\in U^{k}$, $k=1,2,3$,
\begin{gather}
Y\left(u^{2},z\right)u^{1}=\mathcal{I}_{2,1}^{2}(u^{2},z)u^{1};\nonumber \\
Y\left(u^{3},z\right)u^{1}=\mathcal{I}_{3,1}^{3}\left(u^{3},z\right)u^{1};\nonumber \\
Y\left(u^{2},z\right)v^{2}=\left(\mathcal{I}_{2,2}^{1}\left(u^{2},z\right)+\mathcal{I}_{2,2}^{3}\left(u^{2},z\right)\right)v^{2};\nonumber \\
Y\left(u^{2},z\right)v^{3}=\left(\mathcal{I}_{2,3}^{2}\left(u^{2},z\right)+\mathcal{I}_{2,3}^{3}\left(u^{2},z\right)\right)v^{3};\nonumber \\
Y\left(v^{3},z\right)u^{2}=\left(\mathcal{I}_{3,2}^{2}\left(v^{3},z\right)+\mathcal{I}_{3,2}^{3}\left(v^{3},z\right)\right)u^{2};\nonumber \\
Y\left(u^{3},z\right)v^{3}=\left(\mathcal{I}_{3,3}^{1}\left(u^{3},z\right)+\mathcal{I}_{3,3}^{2}\left(u^{3},z\right)+\mathcal{I}_{3,3}^{3}\left(u^{3},z\right)\right)v^{3}\label{U,Y}
\end{gather}

The following result will be applied to prove the uniqueness of the
vertex operator algebra structure on $\mathcal{U}$.

\begin{lemma}\label{new structure} Let $(V,Y)$ be a vertex operator
algebra and $f:V\to V$ be a linear isomorphism such that $\sigma\left(1\right)=1,\sigma\left(\omega\right)=\omega$.
Then $(V,Y^{\sigma})$ is a vertex operator where
\[
Y^{\sigma}(u,z)=\sigma Y(\sigma^{-1}u,z)\sigma^{-1}
\]
and $(V,Y)\cong(V,Y^{\sigma})$. \end{lemma}
\begin{proof}
1) Truncation property: For any $u,v\in V$,
\[
Y^{\sigma}(u,z)=\sigma Y\left(\sigma^{-1}u,z\right)\sigma^{-1}v=\sum_{m\in\mathbb{Z}}\sigma\left(\sigma^{-1}u\right)_{m}\left(\sigma^{-1}v\right)z^{-m-1}.
\]
By the truncation property of $Y$, we have $\left(\sigma^{-1}u\right)_{m}\left(\sigma^{-1}v\right)=0$
for $m\gg0$. Thus $Y^{\sigma}$ satisfies truncation property.

2) Vacuum property:
\[
\lim_{z\to0}Y^{\sigma}\left(u,z\right)1=\lim_{z\to0}\sigma Y\left(\sigma^{-1}u,z\right)\sigma1=\lim_{z\to0}\sigma Y\left(\sigma^{-1}u,z\right)1=\sigma\cdot\sigma^{-1}u=u.
\]

3) $L^{\sigma}(-1)$-derivation property: For any $u\in\mathcal{U}$,
\begin{alignat*}{1}
 & \left[L^{\sigma}\left(-1\right),Y^{\sigma}\left(u,z\right)\right]\\
 & =\left[\sigma L\left(-1\right)\sigma^{-1},\sigma Y\left(\sigma^{-1}u,z\right)\sigma^{-1}\right]\\
 & =\sigma\left[L\left(-1\right),Y\left(\sigma^{-1}u,z\right)\right]\sigma^{-1}\\
 & =\frac{d}{dz}\sigma Y\left(\sigma^{-1}u,z\right)\sigma^{-1}\\
 & =\frac{d}{dz}Y^{\sigma}\left(u,z\right).
\end{alignat*}

4) Commutativity: For any $u,v\in V$,
\[
\left(z_{1}-z_{2}\right)^{m}\left[Y^{\sigma}\left(u,z_{1}\right),Y^{\sigma}\left(v,z_{2}\right)\right]=\left(z_{1}-z_{2}\right)^{m}\sigma\left[Y\left(\sigma^{-1}u,z_{1}\right),Y\left(\sigma^{-1}v,z_{2}\right)\right]\sigma^{-1}=0.
\]

Thus $\left(V,Y^{\sigma}\right)$ is a vertex operator algebra. Since
$\sigma Y\left(u,z\right)\sigma^{-1}=\sigma Y\left(\sigma^{-1}\sigma u,z\right)\sigma^{-1}=Y^{\sigma}(\sigma u,z)$,
we get $\left(V,Y^{\sigma}\right)\cong\left(V,Y\right)$.
\end{proof}
\begin{theorem}The vertex operator algebra structure on $\mathcal{U}$
over $\mathbb{C}$ is unique.\end{theorem}
\begin{proof}
Let $\left(\mathcal{U},Y\right)$ be the vertex operator algebra structure
as given in (\ref{U,Y}). Suppose $\left(\mathcal{U},\overline{Y}\right)$
is another vertex operator algebra structure on $\mathcal{U}$. Without
loss of generality, we may assume $Y\left(u,z\right)=\overline{Y}\left(u,z\right)$
for all $u\in U^{1}$. From our settings above, there exist nonzero
constants $\lambda_{i,1}^{i}$, $\lambda_{2,2}^{j}$, $\lambda_{2,3}^{k}$,
$\lambda_{3,2}^{p}$, $\lambda_{33}^{l}$ where $i,k,p=2,3$, $j=1,3$,
$l=1,2,3$ such that for any $u^{i},v^{i}\in U^{i}$, $i=1,2,3$,
we have
\[
\overline{Y}\left(u^{2},z\right)u^{1}=\lambda_{2,1}^{2}\cdot\mathcal{I}_{2,1}^{2}(u^{2},z)u^{1},
\]
\[
\overline{Y}\left(u^{3},z\right)u^{1}=\lambda_{3,1}^{3}\cdot\mathcal{I}_{3,1}^{3}\left(u^{3},z\right)u^{1},
\]
\[
\overline{Y}\left(u^{2},z\right)v^{2}=\left(\lambda_{2,2}^{1}\cdot\mathcal{I}_{2,2}^{1}\left(u^{2},z\right)+\lambda_{2,2}^{3}\cdot\mathcal{I}_{2,2}^{3}\left(u^{2},z\right)\right)v^{2},
\]
\[
\overline{Y}\left(u^{3},z\right)v^{3}=\left(\lambda_{3,3}^{1}\cdot\mathcal{I}_{3,3}^{1}\left(u^{3},z\right)+\lambda_{3,3}^{2}\cdot\mathcal{I}_{3,3}^{2}\left(u^{3},z\right)+\lambda_{3,3}^{2}\cdot\mathcal{I}_{3,3}^{3}\left(u^{3},z\right)\right)v^{3},
\]
\[
\overline{Y}\left(u^{2},z\right)u^{3}=\left(\lambda_{2,3}^{2}\cdot\mathcal{I}_{2,3}^{2}\left(u^{2},z\right)+\lambda_{2,3}^{3}\cdot\mathcal{I}_{2,3}^{3}\left(u^{2},z\right)\right)u^{3},
\]
\[
\overline{Y}\left(u^{3},z\right)p=\left(\lambda_{3,2}^{2}\cdot\mathcal{I}_{3,2}^{2}\left(u^{3},z\right)+\lambda_{3,2}^{3}\cdot\mathcal{I}_{3,2}^{3}\left(u^{3},z\right)\right)u^{2},
\]
where $\mathcal{I}_{a,b}^{c}\in I_{U^{1}}\left(_{U^{a}\ U^{b}}^{U^{c}}\right),$
$a,b,c\in\left\{ 1,2,3\right\} $ are nonzero intertwining operators.

\emph{Claim 1) $\lambda_{2,1}^{2}=\lambda_{3,1}^{3}=1$}

For any $u^{1}\in U^{1}$, $u^{2}\in U^{2}$, skew symmetry of $Y\left(\cdot,z\right)$
and $\overline{Y}\left(\cdot,z\right)$ ( \cite{FHL} ) imply

\[
\overline{Y}(u^{2},z)u^{1}=e^{zL\left(-1\right)}\overline{Y}\left(u^{1},-z\right)u^{2}=e^{zL\left(-1\right)}Y(u^{1},-z)u^{2}=Y\left(u^{2},z\right)u^{1}=\mathcal{I}_{2,1}^{2}\left(u^{2},z\right)u^{1}.
\]
In the mean time, $\overline{Y}\left(u^{2},z\right)u^{1}=\lambda_{2,1}^{2}\cdot\mathcal{I}_{2,1}^{2}(u^{2},z)u^{1}$.
Thus we get $\lambda_{2,1}^{2}=1$. Similarly, we can prove $\lambda_{3,1}^{1}=1$.

\emph{Claim 2)} $\lambda_{2,2}^{1}=\lambda_{3,3}^{1}=1$.

Note that by Remark \ref{self-dual}, $\mathcal{U}$ has a unique
invariant bilinear form $\left\langle \cdot,\cdot\right\rangle $
with $\left\langle 1,1\right\rangle =1$. For $u^{1}\in U^{1}$ and
$u^{2},v^{2}\in U^{2}$, we have
\[
\left\langle Y\left(u^{2},z)v^{2}\right),u^{1}\right\rangle =\left\langle v^{2},Y\left(e^{zL\left(-1\right)}\left(-z^{-2}\right)^{L\left(0\right)}u^{2},z^{-1}\right)u^{1}\right\rangle .
\]
That is,
\[
\left\langle \mathcal{I}_{2,2}^{1}\left(u^{2},z\right)v^{2},u^{1}\right\rangle =\left\langle v^{2},\mathcal{I}_{2,1}^{2}\left(e^{zL\left(-1\right)}\left(-z^{-2}\right)^{L\left(0\right)}u^{2},z^{-1}\right)u^{1}\right\rangle .
\]
The invariant bilinear form on $\left(\mathcal{U},\overline{Y}\right)$
gives
\[
\left\langle \lambda_{2,2}^{1}\cdot\mathcal{I}_{2,2}^{1}\left(u^{2},z\right)v^{2},u^{1}\right\rangle =\left\langle v^{2},\lambda_{2,1}^{2}\cdot\mathcal{I}_{2,1}^{2}\left(e^{zL\left(-1\right)}\left(-z^{-2}\right)^{L\left(0\right)}u^{2},z^{-1}\right)u^{1}\right\rangle .
\]
Using claim 1, we get $\lambda_{2,2}^{1}=1$. Similarly, we can prove
$\lambda_{3,3}^{1}=1$.

\emph{Claim 3) $\lambda_{2,3}^{k}=\lambda_{3,2}^{k},k=2,3$. }

Let $u^{2},v^{2}\in U^{2}$, $u^{3}\in U^{3}$, by skew symmetry of
$Y$ we obtain
\[
\left\langle Y\left(u^{2},z\right)u^{3},v^{2}\right\rangle =\left\langle e^{zL\left(-1\right)}Y\left(u^{3},-z\right)u^{2},v^{2}\right\rangle ,
\]
that is,
\[
\left\langle \mathcal{I}_{2,3}^{2}\left(u^{2},z\right)u^{3},v^{2}\right\rangle =\left\langle e^{zL\left(-1\right)}\mathcal{I}_{3,2}^{2}\left(u^{3},-z\right)u^{2},v^{2}\right\rangle .
\]
Skew symmetry of $\overline{Y}$ gives
\[
\lambda_{2,3}^{2}\left\langle \mathcal{I}_{2,3}^{2}\left(u^{2},z\right)u^{3},v^{2}\right\rangle =\lambda_{3,2}^{2}\left\langle e^{zL\left(-1\right)}\mathcal{I}_{3,2}^{2}\left(u^{3},-z\right)u^{2},v^{2}\right\rangle .
\]
Comparing the last two identities, we get $\lambda_{2,3}^{2}=\lambda_{3,2}^{2}.$
Similarly, we can prove $\lambda_{2,3}^{3}=\lambda_{3,2}^{3}.$

Let $u^{1}\in U^{1}$, $u^{2},v^{2}\in U^{2}$ and $u^{3}\in U^{3}$,
commutativity of $Y$ and $\overline{Y}$ in (\ref{communitivity})
gives

\[
\iota_{12}^{-1}\left\langle u^{1},\mathcal{I}_{2,2}^{1}\left(u^{2},z_{1}\right)\mathcal{I}_{3,2}^{2}\left(u^{3},z_{2}\right)v^{2}\right\rangle =\iota_{21}^{-1}\left\langle u^{1},\mathcal{I}_{3,3}^{1}\left(u^{3},z_{2}\right)\mathcal{I}_{2,2}^{3}\left(u^{2},z_{1}\right)v^{2}\right\rangle ,
\]

\[
\iota_{12}^{-1}\left\langle u^{1},\lambda_{2,2}^{1}\lambda_{3,2}^{2}\cdot\mathcal{I}_{2,2}^{1}\left(u^{2},z_{1}\right)\mathcal{I}_{3,2}^{2}\left(u^{3},z_{2}\right)v^{2}\right\rangle =\iota_{21}^{-1}\left\langle u^{1},\lambda_{3,3}^{1}\lambda_{2,2}^{3}\cdot\mathcal{I}_{3,3}^{1}\left(u^{3},z_{2}\right)\mathcal{I}_{2,2}^{3}\left(u^{2},z_{1}\right)v^{2}\right\rangle .
\]
The above two identities and claim 2) together give us
\begin{equation}
\lambda_{3,2}^{2}=\lambda_{2,2}^{3}.\label{VOAstr-1}
\end{equation}
Similarly, when $u^{2},v^{2},w^{2}\in U^{2}$, $u^{3}\in U^{3}$,
(\ref{communitivity}) gives
\[
\iota_{12}^{-1}\left\langle w^{2},\mathcal{I}_{2,3}^{2}\left(u^{2},z_{1}\right)\mathcal{I}_{3,2}^{3}\left(u^{3},z_{2}\right)v^{2}\right\rangle =\iota_{21}^{-1}\left\langle w^{2},\mathcal{I}_{3,3}^{2}\left(u^{3},z_{2}\right)\mathcal{I}_{2,2}^{3}\left(u^{2},z_{1}\right)v^{2}\right\rangle ,
\]
\[
\iota_{12}^{-1}\left\langle w^{2},\lambda_{2,3}^{2}\lambda_{3,2}^{3}\cdot\mathcal{I}_{23}^{2}\left(u^{2},z_{1}\right)\mathcal{I}_{3,2}^{3}\left(u^{3},z_{2}\right)v^{2}\right\rangle =\iota_{21}^{-1}\left\langle w^{2},\lambda_{3,3}^{2}\lambda_{2,2}^{3}\cdot\mathcal{I}_{3,3}^{2}\left(u^{3},z_{2}\right)\text{\ensuremath{\mathcal{I}}}_{2,2}^{3}\left(u^{2},z_{1}\right)v^{2}\right\rangle .
\]
Hence by claim 3) and (\ref{VOAstr-1}) we get
\[
\lambda_{3,2}^{3}=\lambda_{3,3}^{2}.
\]
Now we have
\[
\lambda_{2,3}^{2}=\lambda_{3,2}^{2}=\lambda_{2,2}^{3}\ \text{\ \ and\ \ }\ \lambda_{2,3}^{3}=\lambda_{3,2}^{3}=\lambda_{3,3}^{2},
\]
which we denote by $\lambda$ and $\mu$ respectively.

\emph{Claim 4) \label{lambda^2-mu^2}} $\lambda^{2}=\mu^{2}$.

Fix a basis $\left\{ \overline{\mathcal{Y}}_{a,b;i}^{c}|i=1,\cdots,N_{Q_{a},Q_{b}}^{Q_{c}}\right\} $
for $I\left(_{Q_{a},Q_{b}}^{Q_{c}}\right)$, $a,b,c\in\left\{ 1,2,3\right\} $
as in \cite{FFK}. Now we consider the four point functions on $\left(U^{2},U^{3},U^{3},U^{2}\right).$
Let $B_{2,2}^{3,3}$ be as defined in (\ref{B}). Let $t_{1}^{2}\otimes t_{2}^{2}$,
$p_{1}^{2}\otimes p_{2}^{2}\in U^{2},$$w_{1}^{3}\otimes w_{2}^{3},$
$u_{1}^{3}\otimes u_{2}^{3}\in U^{3},$ we have
\begin{equation}
\begin{aligned} & E\left\langle t_{1}^{2}\otimes t_{2}^{2},Y\left(w_{1}^{3}\otimes w_{2}^{3},z_{1}\right)Y\left(u_{1}^{3}\otimes u_{2}^{3},z_{2}\right)p_{1}^{2}\otimes p_{2}^{2}\right\rangle \\
 & =E\langle t_{1}^{2}\otimes t_{2}^{2},\mathcal{Y}_{3,2}^{2}\otimes\overline{\mathcal{Y}}_{3,2}^{2}\left(w_{1}^{3}\otimes w_{2}^{3},z_{1}\right)\cdot\mathcal{Y}_{3,2}^{2}\otimes\overline{\mathcal{Y}}_{3,2}^{2}\left(u_{1}^{3}\otimes u_{2}^{3},z_{2}\right)\cdot p_{1}^{2}\otimes p_{2}^{2}\\
 & +\mathcal{Y}_{3,3}^{2}\otimes\overline{\mathcal{Y}}_{3,2}^{3}\left(w_{1}^{3}\otimes w_{2}^{3},z_{1}\right)\cdot\mathcal{Y}_{3,3}^{2}\otimes\overline{\mathcal{Y}}_{3,2}^{3}\left(u_{1}^{3}\otimes u_{2}^{3},z_{2}\right)\cdot p_{1}^{2}\otimes p_{2}^{2}\rangle\\
 & =E\langle t_{1}^{2}\otimes t_{2}^{2},\mathcal{Y}_{3,2}^{2}\left(w_{1}^{3},z_{1}\right)\mathcal{Y}_{3,2}^{2}\left(u_{1}^{3},z_{2}\right)p_{1}^{2}\otimes\overline{\mathcal{Y}}_{3,2}^{2}\left(w_{2}^{3},z_{1}\right)\overline{\mathcal{Y}}_{3,2}^{2}\left(u_{2}^{3},z_{2}\right)p_{2}^{2}\\
 & +\mathcal{Y}_{3,3}^{2}\left(w_{1}^{3},z_{1}\right)\mathcal{Y}_{3,2}^{3}\left(u_{1}^{3},z_{2}\right)p_{1}^{2}\otimes\overline{\mathcal{Y}}_{3,3}^{2}\left(w_{2}^{3},z_{1}\right)\overline{\mathcal{Y}}_{3,2}^{3}\left(u_{2}^{3},z_{2}\right)p_{2}^{2}\rangle\\
 & =E\langle t_{1}^{2}\otimes t_{2}^{2},\sum_{i=2,3}\left(B_{2,2}^{3,3}\right)_{2,i}\mathcal{Y}_{3,i}^{2}\left(u_{1}^{3},z_{2}\right)\mathcal{Y}_{3,2}^{i}\left(w_{1}^{3},z_{1}\right)p_{1}^{2}\\
 & \otimes\sum_{j=2,3}\left(\tilde{B}_{2,2}^{3,3}\right)_{2,j}\mathcal{\overline{Y}}_{3,j}^{2}\left(u_{2}^{3},z_{2}\right)\mathcal{\overline{Y}}_{3,2}^{j}\left(w_{2}^{3},z_{1}\right)p_{2}^{2}\\
 & +\sum_{i=2,3}\left(B_{2,2}^{3,3}\right)_{3,i}\mathcal{Y}_{3,i}^{2}\left(u_{1}^{3},z_{2}\right)\mathcal{Y}_{3,2}^{i}\left(w_{1}^{3},z_{1}\right)p_{1}^{2}\\
 & \otimes\sum_{j=2,3}\left(\tilde{B}_{2,2}^{3,3}\right)_{2,j}\mathcal{\overline{Y}}_{3,j}^{2}\left(u_{2}^{3},z_{2}\right)\mathcal{\overline{Y}}_{3,2}^{j}\left(w_{2}^{3},z_{1}\right)p_{2}^{2}\rangle
\end{aligned}
\label{2332-1}
\end{equation}
In the mean time, we also have
\begin{alignat}{1}
 & E\left\langle t_{1}^{2}\otimes t_{2}^{2},Y\left(u_{1}^{3}\otimes u_{2}^{3},z_{2}\right)Y\left(w_{1}^{3}\otimes w_{2}^{3},z_{1}\right)p_{1}^{2}\otimes p_{2}^{2}\right\rangle \nonumber \\
 & =E\langle t_{1}^{2}\otimes t_{2}^{2},\mathcal{Y}_{3,2}^{2}\left(u_{1}^{3},z_{2}\right)\mathcal{Y}_{3,2}^{2}\left(w_{1}^{3},z_{1}\right)p_{1}^{2}\otimes\overline{\mathcal{Y}}_{3,2}^{2}\left(u_{2}^{3},z_{2}\right)\overline{\mathcal{Y}}_{3,2}^{2}\left(w_{2}^{3},z_{1}\right)p_{2}^{2}\nonumber \\
 & =+\mathcal{Y}_{3,3}^{2}\left(u_{1}^{3},z_{2}\right)\mathcal{Y}_{3,2}^{3}\left(w_{1}^{3},z_{1}\right)p_{1}^{2}\otimes\overline{\mathcal{Y}}_{3,3}^{2}\left(u_{2}^{3},z_{2}\right)\overline{\mathcal{Y}}_{3,2}^{3}\left(w_{2}^{3},z_{1}\right)p_{2}^{2}\rangle\label{2332-2}
\end{alignat}
Commutativity of $\left(\mathcal{U},Y\right)$ and (\ref{2332-1})
and (\ref{2332-2}) together imply the following system of equations:
\begin{equation}
\begin{cases}
\left(B_{2,2}^{3,3}\right)_{2,2}\cdot\left(\tilde{B}_{2,2}^{3,3}\right)_{2,2}+\left(B_{2,2}^{3,3}\right)_{3,2}\cdot\left(\tilde{B}_{2,2}^{3,3}\right)_{3,2}=1\\
\left(B_{2,2}^{3,3}\right)_{2,3}\cdot\left(\tilde{B}_{2,2}^{3,3}\right)_{2,3}+\left(B_{2,2}^{3,3}\right)_{3,3}\cdot\left(\tilde{B}_{2,2}^{3,3}\right)_{3,3}=1\\
\left(B_{2,2}^{3,3}\right)_{2,2}\cdot\left(\tilde{B}_{2,2}^{3,3}\right)_{2,3}+\left(B_{2,2}^{3,3}\right)_{3,2}\cdot\left(\tilde{B}_{2,2}^{3,3}\right)_{3,3}=0\\
\left(B_{2,2}^{3,3}\right)_{2,3}\cdot\left(\tilde{B}_{2,2}^{3,3}\right)_{2,2}+\left(B_{2,2}^{3,3}\right)_{3,3}\cdot\left(\tilde{B}_{2,2}^{3,3}\right)_{3,2}=0
\end{cases}\label{system (2332)-1}
\end{equation}
Similarly, from commutativity of $\left(\mathcal{U},\overline{Y}\right)$
we have
\begin{align}
 & E\left\langle t_{1}^{2}\otimes t_{2}^{2},\overline{Y}\left(w_{1}^{3}\otimes w_{2}^{3},z_{1}\right)\overline{Y}\left(u_{1}^{3}\otimes u_{2}^{3},z_{2}\right)p_{1}^{2}\otimes p_{2}^{2}\right\rangle \nonumber \\
 & =E\langle t_{1}^{2}\otimes t_{2}^{2},\lambda^{2}\mathcal{Y}_{3,2}^{2}\otimes\overline{\mathcal{Y}}_{3,2}^{2}\left(w_{1}^{3}\otimes w_{2}^{3},z_{1}\right)\cdot\mathcal{Y}_{3,2}^{2}\otimes\overline{\mathcal{Y}}_{3,2}^{2}\left(u_{1}^{3}\otimes u_{2}^{3},z_{2}\right)\cdot p_{1}^{2}\otimes p_{2}^{2}\nonumber \\
 & +\mu^{2}\mathcal{Y}_{3,3}^{2}\otimes\overline{\mathcal{Y}}_{3,2}^{3}\left(w_{1}^{3}\otimes w_{2}^{3},z_{1}\right)\cdot\mathcal{Y}_{3,3}^{2}\otimes\overline{\mathcal{Y}}_{3,2}^{3}\left(u_{1}^{3}\otimes u_{2}^{3},z_{2}\right)\cdot p_{1}^{2}\otimes p_{2}^{2}\rangle\nonumber \\
 & =E\langle t_{1}^{2}\otimes t_{2}^{2},\lambda^{2}\cdot\mathcal{Y}_{3,2}^{2}\left(w_{1}^{3},z_{1}\right)\mathcal{Y}_{3,2}^{2}\left(u_{1}^{3},z_{2}\right)p_{1}^{2}\otimes\overline{\mathcal{Y}}_{3,2}^{2}\left(w_{2}^{3},z_{1}\right)\overline{\mathcal{Y}}_{3,2}^{2}\left(u_{2}^{3},z_{2}\right)p_{2}^{2}\nonumber \\
 & +\mu^{2}\cdot\mathcal{Y}_{3,3}^{2}\left(w_{1}^{3},z_{1}\right)\mathcal{Y}_{3,2}^{3}\left(u_{1}^{3},z_{2}\right)p_{1}^{2}\otimes\overline{\mathcal{Y}}_{3,3}^{2}\left(w_{2}^{3},z_{1}\right)\overline{\mathcal{Y}}_{3,2}^{3}\left(u_{2}^{3},z_{2}\right)p_{2}^{2}\rangle\nonumber \\
 & =E\langle t_{1}^{2}\otimes t_{2}^{2},\lambda^{2}\cdot\sum_{i=2,3}\left(B_{2,2}^{3,3}\right)_{2,i}\mathcal{Y}_{3,i}^{2}\left(u_{1}^{3},z_{2}\right)\mathcal{Y}_{3,2}^{i}\left(w_{1}^{3},z_{1}\right)p_{1}^{2}\nonumber \\
 & \otimes\sum_{j=2,3}\left(\tilde{B}_{2,2}^{3,3}\right)_{2,j}\mathcal{\overline{Y}}_{3,j}^{2}\left(u_{2}^{3},z_{2}\right)\mathcal{\overline{Y}}_{3,2}^{j}\left(w_{2}^{3},z_{1}\right)p_{2}^{2}\nonumber \\
 & +\mu^{2}\cdot\sum_{i=2,3}\left(B_{2,2}^{3,3}\right)_{3,i}\mathcal{Y}_{3,i}^{2}\left(u_{1}^{3},z_{2}\right)\mathcal{Y}_{3,2}^{i}\left(w_{1}^{3},z_{1}\right)p_{1}^{2}\nonumber \\
 & \otimes\sum_{j=2,3}\left(\tilde{B}_{2,2}^{3,3}\right)_{2,j}\mathcal{\overline{Y}}_{3,j}^{2}\left(u_{2}^{3},z_{2}\right)\mathcal{\overline{Y}}_{3,2}^{j}\left(w_{2}^{3},z_{1}\right)p_{2}^{2}\rangle\label{new 2332-1}
\end{align}
and
\begin{alignat}{1}
 & E\left\langle t_{1}^{2}\otimes t_{2}^{2},\overline{Y}\left(u_{1}^{3}\otimes u_{2}^{3},z_{2}\right)\overline{Y}\left(w_{1}^{3}\otimes w_{2}^{3},z_{1}\right)p_{1}^{2}\otimes p_{2}^{2}\right\rangle \nonumber \\
 & =E\langle t_{1}^{2}\otimes t_{2}^{2},\lambda^{2}\cdot\mathcal{Y}_{3,2}^{2}\left(u_{1}^{3},z_{2}\right)\mathcal{Y}_{3,2}^{2}\left(w_{1}^{3},z_{1}\right)p_{1}^{2}\otimes\overline{\mathcal{Y}}_{3,2}^{2}\left(u_{2}^{3},z_{2}\right)\overline{\mathcal{Y}}_{3,2}^{2}\left(w_{2}^{3},z_{1}\right)p_{2}^{2}\nonumber \\
 & +\mu^{2}\cdot\mathcal{Y}_{3,3}^{2}\left(u_{1}^{3},z_{2}\right)\mathcal{Y}_{3,2}^{3}\left(w_{1}^{3},z_{1}\right)p_{1}^{2}\otimes\overline{\mathcal{Y}}_{3,3}^{2}\left(u_{2}^{3},z_{2}\right)\overline{\mathcal{Y}}_{3,2}^{3}\left(w_{2}^{3},z_{1}\right)p_{2}^{2}\rangle\label{new 2332-2}
\end{alignat}
(\ref{new 2332-1}) and (\ref{new 2332-2}) together imply
\begin{equation}
\begin{cases}
\lambda^{2}\cdot\left(B_{2,2}^{3,3}\right)_{2,2}\cdot\left(\tilde{B}_{2,2}^{3,3}\right)_{2,2}+\mu^{2}\cdot\left(B_{2,2}^{3,3}\right)_{3,2}\cdot\left(\tilde{B}_{2,2}^{3,3}\right)_{3,2}=\lambda^{2}\\
\lambda^{2}\cdot\left(B_{2,2}^{3,3}\right)_{2,3}\cdot\left(\tilde{B}_{2,2}^{3,3}\right)_{2,3}+\mu^{2}\cdot\left(B_{2,2}^{3,3}\right)_{3,3}\cdot\left(\tilde{B}_{2,2}^{3,3}\right)_{3,3}=\mu^{2}\\
\lambda^{2}\cdot\left(B_{2,2}^{3,3}\right)_{2,2}\cdot\left(\tilde{B}_{2,2}^{3,3}\right)_{2,3}+\mu^{2}\cdot\left(B_{2,2}^{3,3}\right)_{3,2}\cdot\left(\tilde{B}_{2,2}^{3,3}\right)_{3,3}=0\\
\lambda^{2}\cdot\left(B_{2,2}^{3,3}\right)_{2,3}\cdot\left(\tilde{B}_{2,2}^{3,3}\right)_{2,2}+\mu^{2}\cdot\left(B_{2,2}^{3,3}\right)_{3,3}\cdot\left(\tilde{B}_{2,2}^{3,3}\right)_{3,2}=0
\end{cases}\label{system (2332)-2}
\end{equation}
Systems (\ref{system (2332)-1}) and (\ref{system (2332)-2}) together
imply
\begin{equation}
\begin{cases}
\left(1-\frac{\mu^{2}}{\lambda^{2}}\right)\cdot\left(B_{2,2}^{3,3}\right)_{3,2}\cdot\left(\tilde{B}_{2,2}^{3,3}\right)_{3,2}=0\\
\left(1-\frac{\lambda^{2}}{\mu^{2}}\right)\cdot\left(B_{2,2}^{3,3}\right)_{2,3}\cdot\left(\tilde{B}_{2,2}^{3,3}\right)_{2,3}=0\\
\left(1-\frac{\mu^{2}}{\lambda^{2}}\right)\cdot\left(B_{2,2}^{3,3}\right)_{3,2}\cdot\left(\tilde{B}_{2,2}^{3,3}\right)_{3,3}=0\\
\left(1-\frac{\mu^{2}}{\lambda^{2}}\right)\cdot\left(B_{2,2}^{3,3}\right)_{3,3}\cdot\left(\tilde{B}_{2,2}^{3,3}\right)_{3,2}=0
\end{cases}\label{(2332)-results-1}
\end{equation}
Note that from Lemma \ref{(2332)_3,2 nonzero}, $\left(\tilde{B}_{2,2}^{3,3}\right)_{3,2}\not=0$.
If $\lambda^{2}\not=\mu^{2},$ then $\left(B_{2,2}^{3,3}\right)_{3,2}=\left(B_{2,2}^{3,3}\right)_{3,3}=0$
by the first and fourth identity of (\ref{(2332)-results-1}). Combining
the first identity in (\ref{system (2332)-1}), we obtain $\left(B_{2,2}^{3,3}\right)_{2,2}\not=0$
and $\left(\tilde{B}_{2,2}^{3,3}\right)_{2,2}\not=0$. Combining $\left(B_{2,2}^{3,3}\right)_{2,2}\not=0$,
$\left(B_{2,2}^{3,3}\right)_{3,2}=0$ and the third equality of (\ref{system (2332)-1}),
we get $\left(\tilde{B}_{2,2}^{3,3}\right)_{2,3}=0$. But $\left(B_{2,2}^{3,3}\right)_{3,3}=0$
and the second equality of (\ref{system (2332)-1}) together imply
that $\left(\tilde{B}_{2,2}^{3,3}\right)_{2,3}\not=0$. Contradiction
implies $\lambda^{2}=\mu^{2}.$

\emph{Claim 5) $\lambda_{3,3}^{3}\lambda=1$. }

Consider four point functions on $\left(U^{3},U^{3},U^{2},U^{2}\right)$.
For simplicity, we denote $\lambda_{3,3}^{3}$ by $\gamma.$ Applying
similar argument, we obtain systems
\begin{equation}
\begin{cases}
\left(B_{3,2}^{3,2}\right)_{1,2}\cdot\left(\tilde{B}_{3,2}^{3,2}\right)_{1,2}+\left(B_{3,2}^{3,2}\right)_{3,2}\cdot\left(\tilde{B}_{3,2}^{3,2}\right)_{3,2}=1\\
\left(B_{3,2}^{3,2}\right)_{1,3}\cdot\left(\tilde{B}_{3,2}^{3,2}\right)_{1,3}+\left(B_{3,2}^{3,2}\right)_{3,3}\cdot\left(\tilde{B}_{3,2}^{3,2}\right)_{3,3}=1\\
\left(B_{3,2}^{3,2}\right)_{1,2}\cdot\left(\tilde{B}_{3,2}^{3,2}\right)_{1,3}+\left(B_{3,2}^{3,2}\right)_{3,2}\cdot\left(\tilde{B}_{3,2}^{3,2}\right)_{3,3}=0\\
\left(B_{3,2}^{3,2}\right)_{1,3}\cdot\left(\tilde{B}_{3,2}^{3,2}\right)_{1,2}+\left(B_{3,2}^{3,2}\right)_{3,3}\cdot\left(\tilde{B}_{3,2}^{3,2}\right)_{3,2}=0
\end{cases}\label{system(3322)-1}
\end{equation}
and
\begin{equation}
\begin{cases}
\left(B_{3,2}^{3,2}\right)_{1,2}\cdot\left(\tilde{B}_{3,2}^{3,2}\right)_{1,2}+\gamma\lambda\cdot\left(B_{3,2}^{3,2}\right)_{3,2}\cdot\left(\tilde{B}_{3,2}^{3,2}\right)_{3,2}=\lambda^{2}\\
\left(B_{3,2}^{3,2}\right)_{1,3}\cdot\left(\tilde{B}_{3,2}^{3,2}\right)_{1,3}+\gamma\lambda\cdot\left(B_{3,2}^{3,2}\right)_{3,3}\cdot\left(\tilde{B}_{3,2}^{3,2}\right)_{3,3}=\mu^{2}\\
\left(B_{3,2}^{3,2}\right)_{1,2}\cdot\left(\tilde{B}_{3,2}^{3,2}\right)_{1,3}+\gamma\lambda\cdot\left(B_{3,2}^{3,2}\right)_{3,2}\cdot\left(\tilde{B}_{3,2}^{3,2}\right)_{3,3}=0\\
\left(B_{3,2}^{3,2}\right)_{1,3}\cdot\left(\tilde{B}_{3,2}^{3,2}\right)_{1,2}+\gamma\lambda\cdot\left(B_{3,2}^{3,2}\right)_{3,3}\cdot\left(\tilde{B}_{3,2}^{3,2}\right)_{3,2}=0
\end{cases}\label{system (3322)-2}
\end{equation}
The above two systems together give us
\begin{equation}
\begin{cases}
\left(1-\gamma\lambda\right)\cdot\left(B_{3,2}^{3,2}\right)_{3,2}\cdot\left(\tilde{B}_{3,2}^{3,2}\right)_{3,2}=1-\lambda^{2}\\
\left(1-\gamma\lambda\right)\cdot\left(B_{3,2}^{3,2}\right)_{3,3}\cdot\left(\tilde{B}_{3,2}^{3,2}\right)_{3,3}=1-\mu^{2}\\
\left(1-\gamma\lambda\right)\cdot\left(B_{3,2}^{3,2}\right)_{3,2}\cdot\left(\tilde{B}_{3,2}^{3,2}\right)_{3,3}=0\\
\left(1-\gamma\lambda\right)\cdot\left(B_{3,2}^{3,2}\right)_{3,3}\cdot\left(\tilde{B}_{3,2}^{3,2}\right)_{3,2}=0
\end{cases}\label{*}
\end{equation}
Set
\[
S=\left(\begin{array}{cc}
\left(B_{3,2}^{3,2}\right){}_{1,2} & \left(B_{3,2}^{3,2}\right)_{1,3}\\
\left(B_{3,2}^{3,2}\right)_{3,2} & \left(B_{3,2}^{3,2}\right)_{3,3}
\end{array}\right),T=\left(\begin{array}{cc}
\left(\tilde{B}_{3,2}^{3,2}\right)_{1,2} & \left(\tilde{B}_{3,2}^{3,2}\right)_{3,2}\\
\left(\tilde{B}_{3,2}^{3,2}\right)_{1,3} & \left(\tilde{B}_{3,2}^{3,2}\right)_{3,3}
\end{array}\right).
\]
Then system (\ref{system(3322)-1}) implies $S^{T}T=\left(\begin{array}{cc}
1 & 0\\
0 & 1
\end{array}\right)$. So $T^{-1}=S^{T}$, which gives
\begin{equation}
\frac{1}{\det T}\left(\begin{array}{cc}
\left(\tilde{B}_{3,2}^{3,2}\right)_{3,3} & -\left(\tilde{B}_{3,2}^{3,2}\right)_{3,2}\\
-\left(\tilde{B}_{3,2}^{3,2}\right)_{1,3} & \left(\tilde{B}_{3,2}^{3,2}\right)_{1,2}
\end{array}\right)=\left(\begin{array}{cc}
\left(B_{3,2}^{3,2}\right)_{1,2} & \left(B_{3,2}^{3,2}\right)_{3,2}\\
\left(B_{3,2}^{3,2}\right)_{1,3} & \left(B_{3,2}^{3,2}\right)_{3,3}
\end{array}\right).\label{IV(3322)}
\end{equation}
From Lemma \ref{(3322)_3,3 nonzero} $\left(\tilde{B}_{3,2}^{3,2}\right)_{3,3}\not=0$.
Using (\ref{IV(3322)}), we get $\left(B_{3,2}^{3,2}\right)_{1,2}\not=0$.
Assume that $\gamma\lambda\not=1$, then $\left(\tilde{B}_{3,2}^{3,2}\right)_{3,3}\not=0$
and the third equation in (\ref{*}) together imply that $\left(B_{3,2}^{3,2}\right)_{3,2}=0$.
So we have $\left(B_{3,2}^{3,2}\right)_{3,2}\cdot\left(\tilde{B}_{3,2}^{3,2}\right)_{3,2}=0.$
Note that $\left(B_{3,2}^{3,2}\right)_{1,2}\not=0,$ $\left(B_{3,2}^{3,2}\right)_{3,2}=0$
and the third equation in (\ref{system(3322)-1}) together imply that
$\left(\tilde{B}_{3,2}^{3,2}\right)_{1,3}=0.$ By second equation
in (\ref{system(3322)-1}), we have $\left(B_{3,2}^{3,2}\right)_{3,3}\cdot\left(\tilde{B}_{3,2}^{3,2}\right)_{3,3}=1.$
Since we have proved that $\lambda^{2}=\mu^{2}$, the first two equations
of System (\ref{*}) imply that $\left(B_{3,2}^{3,2}\right)_{3,2}\cdot\left(\tilde{B}_{3,2}^{3,2}\right)_{3,2}=\left(B_{3,2}^{3,2}\right)_{3,3}\cdot\left(\tilde{B}_{3,2}^{3,2}\right)_{3,3}$
. Contradiction implies $\gamma\lambda=1$.

\emph{Claim 6) $\lambda=\gamma$. }

Consider four point functions on $\left(U^{2},U^{3},U^{3},U^{3}\right)$.
Apply similar arguments as above on $\left(\mathcal{U},Y\right)$
and $\left(\mathcal{U},\overline{Y}\right)$, we obtain the following
systems respectively:
\begin{equation}
\begin{cases}
\left(B_{2,3}^{3,3}\right)_{2,2}\cdot\left(\tilde{B}_{2,3}^{3,3}\right)_{2,2}+\left(B_{2,3}^{3,3}\right)_{3,2}\cdot\left(\tilde{B}_{2,3}^{3,3}\right)_{3,2}=1\\
\left(B_{2,3}^{3,3}\right)_{2,3}\cdot\left(\tilde{B}_{2,3}^{3,3}\right)_{2,3}+\left(B_{2,3}^{3,3}\right)_{3,3}\cdot\left(\tilde{B}_{2,3}^{3,3}\right)_{3,3}=1\\
\left(B_{2,3}^{3,3}\right)_{2,2}\cdot\left(\tilde{B}_{2,3}^{3,3}\right)_{2,3}+\left(B_{2,3}^{3,3}\right)_{3,2}\cdot\left(\tilde{B}_{2,3}^{3,3}\right)_{3,3}=0\\
\left(B_{2,3}^{3,3}\right)_{2,3}\cdot\left(\tilde{B}_{2,3}^{3,3}\right)_{2,2}+\left(B_{2,3}^{3,3}\right)_{3,3}\cdot\left(\tilde{B}_{2,3}^{3,3}\right)_{3,2}=0
\end{cases}\label{system (2333)-1}
\end{equation}
\begin{equation}
\begin{cases}
\lambda\mu\cdot\left(B_{2,3}^{3,3}\right)_{2,2}\cdot\left(\tilde{B}_{2,3}^{3,3}\right)_{2,2}+\mu\gamma\cdot\left(B_{2,3}^{3,3}\right)_{3,2}\cdot\left(\tilde{B}_{2,3}^{3,3}\right)_{3,2}=\lambda\mu\\
\lambda\mu\cdot\left(B_{2,3}^{3,3}\right)_{2,3}\cdot\left(\tilde{B}_{2,3}^{3,3}\right)_{2,3}+\mu\gamma\cdot\left(B_{2,3}^{3,3}\right)_{3,3}\cdot\left(\tilde{B}_{2,3}^{3,3}\right)_{3,3}=\mu\gamma\\
\lambda\mu\cdot\left(B_{2,3}^{3,3}\right)_{2,2}\cdot\left(\tilde{B}_{2,3}^{3,3}\right)_{2,3}+\mu\gamma\cdot\left(B_{2,3}^{3,3}\right)_{3,2}\cdot\left(\tilde{B}_{2,3}^{3,3}\right)_{3,3}=0\\
\lambda\mu\cdot\left(B_{2,3}^{3,3}\right)_{2,3}\cdot\left(\tilde{B}_{2,3}^{3,3}\right)_{2,2}+\mu\gamma\cdot\left(B_{2,3}^{3,3}\right)_{3,3}\cdot\left(\tilde{B}_{2,3}^{3,3}\right)_{3,2}=0
\end{cases}\label{system (2333)-2}
\end{equation}
(\ref{system (2333)-1}) and (\ref{system (2333)-2}) together gives
\begin{equation}
\begin{cases}
\left(1-\frac{\gamma}{\lambda}\right)\cdot\left(B_{2,3}^{3,3}\right)_{3,2}\cdot\left(\tilde{B}_{2,3}^{3,3}\right)_{3,2}=0\\
\left(1-\frac{\lambda}{\gamma}\right)\cdot\left(B_{2,3}^{3,3}\right)_{2,3}\cdot\left(\tilde{B}_{2,3}^{3,3}\right)_{2,3}=0\\
\left(1-\frac{\gamma}{\lambda}\right)\cdot\left(B_{2,3}^{3,3}\right)_{3,2}\cdot\left(\tilde{B}_{2,3}^{3,3}\right)_{3,3}=0\\
\left(1-\frac{\gamma}{\lambda}\right)\cdot\left(B_{2,3}^{3,3}\right)_{3,3}\cdot\left(\tilde{B}_{2,3}^{3,3}\right)_{3,2}=0
\end{cases}\label{(2333) result}
\end{equation}
By Lemma \ref{(2333)_3,2 nonzero} $\left(\tilde{B}_{2,3}^{3,3}\right)_{3,2}\not=0$.
Assume $\frac{\gamma}{\lambda}\not=1$, then the fourth equation of
(\ref{(2333) result}) imply $\left(B_{2,3}^{3,3}\right)_{3,3}=0$.
Using the second equations in (\ref{system (2333)-1}) we get $\left(B_{2,3}^{3,3}\right)_{2,3}\cdot\left(\tilde{B}_{2,3}^{3,3}\right)_{2,3}=1$,
which contradicts with the second equation of (\ref{(2333) result}).
Therefore, $\frac{\gamma}{\lambda}=1$.

The above claims together imply
\[
\lambda=\mu=\gamma=\pm1
\]
or
\[
\lambda=\gamma=1,\mu=-1
\]
or
\[
\lambda=\gamma=-1,\mu=1.
\]
Define a linear map $\sigma$ such that
\[
\sigma|_{U^{1}}=1,\ \sigma|_{U^{2}}=\mu,\ \sigma|_{U^{3}}=\lambda
\]
where $\lambda=\pm1$ and $\mu=\pm1$. It is clear that $\sigma$
is a linear isomorphism of $\mathcal{U}$. Using Lemma \ref{new structure},
$\sigma$ gives a vertex operator algebra structure $(\mathcal{U},Y^{\sigma})$
with $Y^{\sigma}(u,z)=\sigma Y(\sigma^{-1}u,z)\sigma^{-1}$ which
is isomorphic to $(\mathcal{U},Y)$. It is easy to verify that $Y^{\sigma}\left(u,z\right)=\overline{Y}\left(u,z\right)$
for all $u\in\mathcal{U}$. Thus we proved the uniqueness of the vertex
operator algebra structure on $\mathcal{U}$ .
\end{proof}

\section{Classification of irreducible modules}

In this section, we will classify all the irreducible modules for
$\mathcal{U}$. First we will find 14 irreducible $\mathcal{U}$-modules.
To show they give all the irreducible modules, we shall use the theory
of quantum dimensions. For simplicity, we shall use $[h_{1},\ h_{2}]$
to denote the module $\mathcal{V}(h_{1})\otimes L(\frac{25}{28},\ h_{2})$.

\subsection{Realization of irreducible $\mathcal{U}$-modules}

Let $A_{1}=\mathbb{Z}\alpha$, with $\left\langle \alpha,\alpha\right\rangle =2$,
be the root lattice of type $A_{1}$ and $V_{A_{1}}$ the lattice
vertex operator algebra associated with $A_{1}$. It is well known
that the irreducible $V_{A_{1}}$-modules $V_{A_{1}}$ and $V_{\frac{\alpha}{2}+A_{1}}$
are both level one representations of $\hat{\mathfrak{sl}_{2}}(\mathbb{C})$
\cite{DL,FLM}. In fact, $V_{A_{1}}\cong\mathcal{L}\left(1,0\right)$
and $V_{\frac{\alpha}{2}+A_{1}}\cong\mathcal{L}\left(1,1\right)$.
Let $V_{A_{1}^{6}}$ be the lattice vertex operator algebra associated
with the lattice $A_{1}^{6}$, where $A_{1}^{6}=\mathbb{Z}\alpha_{1}\oplus\cdots\oplus\mathbb{Z}\alpha_{6}$
is the orthogonal sum of 6 copies of $A_{1}$. Then we have
\[
V_{A_{1}^{6}}\cong V_{A_{1}}\otimes\cdots\otimes V_{A_{1}}\cong\mathcal{L}\left(1,0\right){}^{\otimes6}
\]
as a vertex operator algebra and
\[
V_{\gamma+A_{1}^{6}}\cong\mathcal{L}\left(1,1\right)^{\otimes4}\otimes\mathcal{L}\left(1,0\right)\otimes\mathcal{L}\left(1,0\right)
\]
as a module of $\mathcal{L}(1,\ 0)^{\otimes6}$ , where $\gamma=\frac{1}{2}\alpha_{1}+\frac{1}{2}\alpha_{2}+\frac{1}{2}\alpha_{3}+\frac{1}{2}\alpha_{4}$.
Set $L=A_{1}^{6}\cup(\gamma+A_{1}^{6})$, then $L$ is an even lattice
and we have an isomorphism
\[
V_{L}=V_{A_{1}^{6}}\oplus V_{\gamma+A_{1}^{6}}\cong\left\{ \mathcal{L}\left(1,0\right)^{\otimes4}\oplus\mathcal{L}\left(1,1\right)^{\otimes4}\right\} \otimes\mathcal{L}\left(1,0\right)\otimes\mathcal{L}\left(1,0\right).
\]
Using (\ref{eq:(2.3)}) we have the following inclusions
\[
\mathcal{L}\left(1,0\right)^{\otimes3}\supset L\left(\frac{1}{2},0\right)\otimes L\left(\frac{7}{10},0\right)\otimes\mathcal{L}\left(3,0\right),
\]
\[
\mathcal{L}\left(1,1\right)^{\otimes3}\supset L\left(\frac{1}{2},0\right)\otimes L\left(\frac{7}{10},0\right)\otimes\mathcal{L}\left(3,3\right).
\]
Thus, $V_{L}$ contains a vertex operator subalgebra isomorphic to
\[
\mathcal{L}\left(3,0\right)\otimes\mathcal{L}\left(1,0\right)\otimes\mathcal{L}\left(1,0\right)\otimes\mathcal{L}\left(1,0\right)\oplus\mathcal{L}\left(3,3\right)\otimes\mathcal{L}\left(1,1\right)\otimes\mathcal{L}\left(1,0\right)\otimes\mathcal{L}\left(1,0\right).
\]

By (\ref{eq:(2.3)}) and straightforward calculation, we get the following
lemma:

\begin{lemma} \label{realization of VOA}We have the following decomposition:
\begin{align*}
 & \mathcal{L}\left(3,0\right)\otimes\mathcal{L}\left(1,0\right)\otimes\mathcal{L}\left(1,0\right)\otimes\mathcal{L}\left(1,0\right)\oplus\mathcal{L}\left(3,3\right)\otimes\mathcal{L}\left(1,1\right)\otimes\mathcal{L}\left(1,0\right)\otimes\mathcal{L}\left(1,0\right)\\
 & \cong\left\{ \left[0,\ 0\right]\oplus\left[\frac{1}{7},\frac{34}{7}\right]\oplus\left[\frac{5}{7},\frac{9}{7}\right]\right\} \otimes\mathcal{L}\left(6,0\right)\\
 & \oplus\left\{ \left[0,\frac{3}{4}\right]\oplus\left[\frac{5}{7},\frac{1}{28}\right]\oplus\left[\frac{1}{7},\frac{45}{28}\right]\right\} \otimes\mathcal{L}\left(6,2\right)\\
 & \oplus\left\{ \left[0,\frac{13}{4}\right]\oplus\left[\frac{1}{7},\frac{3}{28}\right]\oplus\left[\frac{5}{7},\frac{15}{28}\right]\right\} \otimes\mathcal{L}\left(6,4\right)\\
 & \oplus\left\{ \left[0,\frac{15}{2}\right]\oplus\left[\frac{1}{7},\frac{5}{14}\right]\oplus\left[\frac{5}{7},\frac{39}{14}\right]\right\} \otimes\mathcal{L}\left(6,6\right)
\end{align*}
Thus $\text{\ensuremath{\mathcal{L}}}\left(3,0\right)\otimes\mathcal{L}\left(1,0\right)\otimes\mathcal{L}\left(1,0\right)\otimes\mathcal{L}\left(1,0\right)\oplus\mathcal{L}\left(3,3\right)\otimes\mathcal{L}\left(1,1\right)\otimes\mathcal{L}\left(1,0\right)\otimes\mathcal{L}\left(1,0\right)$
and $V_{L}$ contain a vertex operator subalgebra isomorphic to
\[
\left[0,\ 0\right]\oplus\left[\frac{1}{7},\frac{34}{7}\right]\oplus\left[\frac{5}{7},\frac{9}{7}\right]
\]
which is isomorphic to $\mathcal{U}$ from the uniqueness of $\mathcal{U}$
discussed in Section 3.

\end{lemma}

\begin{lemma} \label{Real modules} The following list give 14 irreducible
$\mathcal{U}$-module.
\begin{align*}
M^{0} & =\left[0,0\right]\oplus\left[\frac{1}{7},\frac{34}{7}\right]\oplus\left[\frac{5}{7},\frac{9}{7}\right], & M^{1} & =\left[0,\frac{3}{4}\right]\oplus\left[\frac{1}{7},\frac{45}{28}\right]\oplus\left[\frac{5}{7},\frac{1}{28}\right],\\
M^{2} & =\left[0,\frac{13}{4}\right]\oplus\left[\frac{1}{7},\frac{3}{28}\right]\oplus\left[\frac{5}{7},\frac{15}{28}\right], & M^{3} & =\left[0,\frac{15}{2}\right]\oplus\left[\frac{1}{7},\frac{5}{14}\right]\oplus\left[\frac{5}{7},\frac{39}{14}\right],\\
M^{4} & =\left[0,\frac{165}{32}\right]\oplus\left[\frac{1}{7},\frac{3}{224}\right]\oplus\left[\frac{5}{7},\frac{323}{224}\right], & M^{5} & =\left[0,\frac{5}{32}\right]\oplus\left[\frac{1}{7},\frac{675}{224}\right]\oplus\left[\frac{5}{7},\frac{99}{224}\right],\\
M^{6} & =\left[0,\frac{57}{32}\right]\oplus\left[\frac{1}{7},\frac{143}{224}\right]\oplus\left[\frac{5}{7},\frac{15}{224}\right], & M^{7} & =\left[\frac{2}{5},0\right]\oplus\left[\frac{19}{35},\frac{34}{7}\right]\oplus\left[\frac{39}{35},\frac{9}{7}\right],\\
M^{8} & =\left[\frac{2}{5},\frac{3}{4}\right]\oplus\left[\frac{19}{35},\frac{45}{28}\right]\oplus\left[\frac{39}{35},\frac{1}{28}\right], & M^{9} & =\left[\frac{2}{5},\frac{13}{4}\right]\oplus\left[\frac{19}{35},\frac{3}{28}\right]\oplus\left[\frac{39}{35},\frac{15}{28}\right],\\
M^{10} & =\left[\frac{2}{5},\frac{15}{2}\right]\oplus\left[\frac{19}{35},\frac{5}{14}\right]\oplus\left[\frac{39}{35},\frac{39}{14}\right], & M^{11} & =\left[\frac{2}{5},\frac{5}{32}\right]\oplus\left[\frac{19}{35},\frac{675}{224}\right]\oplus\left[\frac{39}{35},\frac{99}{224}\right],\\
M^{12} & =\left[\frac{2}{5},\frac{57}{32}\right]\oplus\left[\frac{19}{35},\frac{143}{224}\right]\oplus\left[\frac{39}{35},\frac{15}{224}\right], & M^{13} & =\left[\frac{2}{5},\frac{165}{32}\right]\oplus\left[\frac{19}{35},\frac{3}{224}\right]\oplus\left[\frac{39}{35},\frac{323}{224}\right].
\end{align*}
\end{lemma}
\begin{proof}
From Remark \ref{rationality-1} and Lemma \ref{realization of VOA},
$\mathcal{U}\otimes\mathcal{L}\left(6,0\right)$ is a rational vertex
operator subalgebra of the vertex operator algebra
\[
K=\left\{ \mathcal{L}\left(3,0\right)\otimes\mathcal{L}\left(1,0\right)\oplus\mathcal{L}\left(3,3\right)\otimes\mathcal{L}\left(1,1\right)\right\} \otimes\mathcal{L}\left(1,0\right)\otimes\mathcal{L}\left(1,0\right).
\]
So each irreducible $K$-module is a direct sum of irreducible $\mathcal{U}\otimes\mathcal{L}\left(6,0\right)$-modules.
From Proposition 5.2 \cite{Li1} we know that $\mathcal{L}\left(3,2\right)\otimes\mathcal{L}\left(1,0\right)\oplus\mathcal{L}\left(3,1\right)\otimes\mathcal{L}\left(1,1\right)$
is an irreducible module for $\mathcal{L}\left(3,0\right)\otimes\mathcal{L}\left(1,0\right)\oplus\mathcal{L}\left(3,3\right)\otimes\mathcal{L}\left(1,1\right)$.
Thus we have the following irreducible $K$-modules:
\[
\begin{array}{c}
\left\{ \mathcal{L}\left(3,2\right)\otimes\mathcal{L}\left(1,0\right)\oplus\mathcal{L}\left(3,1\right)\otimes\mathcal{L}\left(1,1\right)\right\} \otimes\mathcal{L}\left(1,0\right)\otimes\mathcal{L}\left(1,0\right),\\
\left\{ \mathcal{L}\left(3,2\right)\otimes\mathcal{L}\left(1,0\right)\oplus\mathcal{L}\left(3,1\right)\otimes\mathcal{L}\left(1,1\right)\right\} \otimes\mathcal{L}\left(1,0\right)\otimes\mathcal{L}\left(1,1\right),\\
\left\{ \mathcal{L}\left(3,0\right)\otimes\mathcal{L}\left(1,0\right)\oplus\mathcal{L}\left(3,3\right)\otimes\mathcal{L}\left(1,1\right)\right\} \otimes\mathcal{L}\left(1,0\right)\otimes\mathcal{L}\left(1,1\right).
\end{array}
\]
Using (\ref{eq:(2.3)}) we obtain the following decomposition:
\begin{alignat*}{1}
 & \mathcal{L}\left(3,0\right)\otimes\mathcal{L}\left(1,0\right)\otimes\mathcal{L}\left(1,0\right)\otimes\mathcal{L}\left(1,0\right)\oplus\mathcal{L}\left(3,3\right)\otimes\mathcal{L}\left(1,1\right)\otimes\mathcal{L}\left(1,0\right)\otimes\mathcal{L}\left(1,0\right)\\
 & \cong\left\{ \left[0,\ 0\right]\oplus\left[\frac{1}{7},\frac{34}{7}\right]\oplus\left[\frac{5}{7},\frac{9}{7}\right]\right\} \otimes\mathcal{L}\left(6,0\right)\\
 & \oplus\left\{ \left[0,\frac{3}{4}\right]\oplus\left[\frac{1}{7},\frac{45}{28}\right]\oplus\left[\frac{5}{7},\frac{1}{28}\right]\right\} \otimes\mathcal{L}\left(6,2\right)\\
 & \oplus\left\{ \left[0,\frac{13}{4}\right]\oplus\left[\frac{1}{7},\frac{3}{28}\right]\oplus\left[\frac{5}{7},\frac{15}{28}\right]\right\} \otimes\mathcal{L}\left(6,4\right)\\
 & \oplus\left\{ \left[0,\frac{15}{2}\right]\oplus\left[\frac{1}{7},\frac{5}{14}\right]\oplus\left[\frac{5}{7},\frac{39}{14}\right]\right\} \otimes\mathcal{L}\left(6,6\right),
\end{alignat*}
\begin{alignat*}{1}
 & \mathcal{L}\left(3,0\right)\otimes\mathcal{L}\left(1,0\right)\otimes\mathcal{L}\left(1,0\right)\otimes\mathcal{L}\left(1,1\right)\oplus\mathcal{L}\left(3,3\right)\otimes\mathcal{L}\left(1,1\right)\otimes\mathcal{L}\left(1,0\right)\otimes\mathcal{L}\left(1,1\right)\\
 & \cong\left\{ \left[0,\frac{5}{32}\right]\oplus\left[\frac{1}{7},\frac{675}{224}\right]\oplus\left[\frac{5}{7},\frac{99}{224}\right]\right\} \otimes\mathcal{L}\left(6,1\right)\\
 & \oplus\left\{ \left[0,\frac{57}{32}\right]\oplus\left[\frac{1}{7},\frac{143}{224}\right]\oplus\left[\frac{5}{7},\frac{15}{224}\right]\right\} \otimes\mathcal{L}\left(6,3\right)\\
 & \oplus\left\{ \left[0,\frac{165}{32}\right]\oplus\left[\frac{1}{7},\frac{3}{224}\right]\oplus\left[\frac{5}{7},\frac{323}{224}\right]\right\} \otimes\mathcal{L}\left(6,5\right),
\end{alignat*}
\begin{alignat*}{1}
 & \mathcal{L}\left(3,2\right)\otimes\mathcal{L}\left(1,0\right)\otimes\mathcal{L}\left(1,0\right)\otimes\mathcal{L}\left(1,1\right)\oplus\mathcal{L}\left(3,1\right)\otimes\mathcal{L}\left(1,1\right)\otimes\mathcal{L}\left(1,0\right)\otimes\mathcal{L}\left(1,1\right)\\
 & \cong\left\{ \left[\frac{2}{5},\frac{5}{32}\right]\oplus\left[\frac{19}{35},\frac{675}{224}\right]\oplus\left[\frac{39}{35},\frac{99}{224}\right]\right\} \otimes\mathcal{L}\left(6,1\right)\\
 & \oplus\left\{ \left[\frac{2}{5},\frac{57}{32}\right]\oplus\left[\frac{19}{35},\frac{143}{224}\right]\oplus\left[\frac{39}{35},\frac{15}{224}\right]\right\} \otimes\mathcal{L}\left(6,3\right)\\
 & \oplus\left\{ \left[\frac{2}{5},\frac{165}{32}\right]\oplus\left[\frac{19}{35},\frac{3}{224}\right]\oplus\left[\frac{39}{35},\frac{323}{224}\right]\right\} \otimes\mathcal{L}\left(6,5\right),
\end{alignat*}
\begin{alignat*}{1}
 & \mathcal{L}\left(3,2\right)\otimes\mathcal{L}\left(1,0\right)\otimes\mathcal{L}\left(1,0\right)\otimes\mathcal{L}\left(1,0\right)\oplus\mathcal{L}\left(3,1\right)\otimes\mathcal{L}\left(1,1\right)\otimes\mathcal{L}\left(1,0\right)\otimes\mathcal{L}\left(1,0\right)\\
 & \cong\left\{ \left[\frac{2}{5},0\right]\oplus\left[\frac{19}{35},\frac{34}{7}\right]\oplus\left[\frac{39}{35},\frac{9}{7}\right]\right\} \otimes\mathcal{L}\left(6,0\right)\\
 & \oplus\left\{ \left[\frac{2}{5},\frac{3}{4}\right]\oplus\left[\frac{19}{35},\frac{45}{28}\right]\oplus\left[\frac{39}{35},\frac{1}{28}\right]\right\} \otimes\mathcal{L}\left(6,2\right)\\
 & \oplus\left\{ \left[\frac{2}{5},\frac{13}{4}\right]\oplus\left[\frac{19}{35},\frac{3}{28}\right]\oplus\left[\frac{39}{35},\frac{15}{28}\right]\right\} \otimes\mathcal{L}\left(6,4\right)\\
 & \oplus\left\{ \left[\frac{2}{5},\frac{15}{2}\right]\oplus\left[\frac{19}{35},\frac{5}{14}\right]\oplus\left[\frac{39}{35},\frac{39}{14}\right]\right\} \otimes\mathcal{L}\left(6,6\right).
\end{alignat*}
Thus we see that $M^{0},M^{1},\cdots,M^{13}$ are $\mathcal{U}$-modules.
It is easy to see that $M^{i},i=0,1,\cdots,13$ are irreducible by
fusion rules of irreducible $L\left(\frac{25}{28},0\right)$-modules
and $\mathcal{V}$-modules in Propositions \ref{fusion rules of virasoro modules}
and \ref{fusion rules for U3A}.
\end{proof}
\begin{remark}\label{module str.,dual} For modules $M^{i},\ i=0,1,\ \cdots,\ 13$
in Lemma \ref{Real modules}, we denote the summands of each $M^{i}$
by $M_{1}^{i},\ M_{2}^{i},\ M_{3}^{i}$ from left to right. Note that
$M_{2}^{i}=U^{2}\boxtimes_{U^{1}}M_{1}^{i}$, $M_{3}^{i}=U^{3}\boxtimes_{U^{1}}M_{1}^{i}$,
$i=0,1,\ \cdots,\ 13$. Thus $M^{i}=\mathcal{U}\boxtimes_{U^{1}}M_{1}^{i}$,
$i=0,1,\ \cdots,\ 13$. Consider quantum dimensions of both sides,
applying Proposition \ref{product property of quantum dimension }
we obtain
\[
q\dim_{U^{1}}M^{i}=q\dim_{U^{1}}\mathcal{U}\cdot q\dim_{U^{1}}M_{1}^{i}
\]
that is, $\frac{q\dim_{U^{1}}M^{i}}{q\dim_{U^{1}}\mathcal{U}}=q\dim_{U^{1}}M_{1}^{i}$
and hence we have
\begin{equation}
q\dim_{\mathcal{U}}M^{i}=q\dim_{U^{1}}M_{1}^{i},i=0,1,\cdots,13.\label{quntum dimension of modules relations}
\end{equation}

\end{remark}

\subsection{Classification}

To finish the classification of irreducible $\mathcal{U}$-modules,
we will show that the list of $\mathcal{U}$-modules in Lemma \ref{Real modules}
give all the irreducible inequivalent $\mathcal{U}$-modules. For
this goal, we will compute global dimension of $\mathcal{U}$.

Using the tables in Section \ref{Fusion rules for P_i's} and properties
of quantum dimensions in Proposition \ref{product property of quantum dimension }
we get
\[
q\dim_{\mathcal{V}}\mathcal{V}\left(\frac{1}{7}\right)=q\dim_{L\left(\frac{25}{28},0\right)}L\left(\frac{25}{28},\frac{34}{7}\right),
\]
\[
q\dim_{\mathcal{V}}\mathcal{V}\left(\frac{5}{7}\right)=q\dim_{L\left(\frac{25}{28},0\right)}L\left(\frac{25}{28},\frac{9}{7}\right).
\]
Denote them by $x$ and $y$ respectively.

From fusion rules in Proposition \ref{fusion rules for U3A}, we see
that $\mathcal{V}\left(\frac{2}{5}\right)\boxtimes_{\mathcal{V}}\mathcal{V}\left(\frac{2}{5}\right)=\mathcal{V}\left(0\right)+\mathcal{V}\left(\frac{2}{5}\right).$
Using property of quantum dimension in Proposition \ref{product property of quantum dimension },
$\left(q\dim_{\mathcal{V}}\mathcal{V}\left(\frac{2}{5}\right)\right)^{2}=1+q\dim_{\mathcal{V}}\mathcal{V}\left(\frac{2}{5}\right)$
and hence we get $q\dim_{\mathcal{V}}\mathcal{V}\left(\frac{2}{5}\right)=\frac{1+\sqrt{5}}{2}$.
Combining the fusion rules for irreducible $\mathcal{V}$-modules
in Propositions \ref{fusion rules for U3A} and \ref{product property of quantum dimension },
one can find
\begin{alignat*}{1}
 & \text{glob}\mathcal{V}\\
 & =1+x^{2}+y^{2}+\left(\frac{1+\sqrt{5}}{2}\right)^{2}+\left(x\cdot\left(\frac{1+\sqrt{5}}{2}\right)\right)^{2}+\left(y\cdot\left(\frac{1+\sqrt{5}}{2}\right)\right)^{2}\\
 & =\left(1+x^{2}+y^{2}\right)\cdot\left(1+\left(\frac{1+\sqrt{5}}{2}\right)^{2}\right).
\end{alignat*}
Recall that the highest weights for irreducible $L\left(\frac{25}{28},0\right)$-modules
are given in Remark \ref{pairs corr. to U's weights for unitary model}
and fusion rules for these irreducible modules are given in Theorem
\ref{fusion rules of virasoro modules}. For convenience, we list
fusion rules explicitly for some irreducible $L\left(\frac{25}{28},0\right)$-modules
which will help us determine quantum dimensions of these modules.
\begin{center}
\begin{tabular}{|c|c|c|}
\hline
$\boxtimes$  & $L\left(\frac{25}{28},\frac{34}{7}\right)$  & $L\left(\frac{25}{28},\frac{9}{7}\right)$\tabularnewline
\hline
\hline
$L\left(\frac{25}{28},\frac{3}{4}\right)$  & $L\left(\frac{25}{28},\frac{45}{28}\right)$  & $L\left(\frac{25}{28},\frac{1}{28}\right)$\tabularnewline
\hline
\multicolumn{1}{|c|}{$L\left(\frac{25}{28},\frac{15}{2}\right)$} & $L\left(\frac{25}{28},\frac{5}{14}\right)$  & $L\left(\frac{25}{28},\frac{39}{14}\right)$\tabularnewline
\hline
$L\left(\frac{25}{28},\frac{165}{32}\right)$  & $L\left(\frac{25}{28},\frac{3}{224}\right)$  & $L\left(\frac{25}{28},\frac{323}{224}\right)$\tabularnewline
\hline
$L\left(\frac{25}{28},\frac{5}{32}\right)$  & $L\left(\frac{25}{28},\frac{675}{224}\right)$  & $L\left(\frac{25}{28},\frac{99}{224}\right)$\tabularnewline
\hline
$L\left(\frac{25}{28},\frac{13}{4}\right)$  & $L\left(\frac{25}{28},\frac{3}{28}\right)$  & $L\left(\frac{25}{28},\frac{15}{28}\right)$\tabularnewline
\hline
$L\left(\frac{25}{28},\frac{57}{32}\right)$  & $L\left(\frac{25}{28},\frac{143}{224}\right)$  & $L\left(\frac{25}{28},\frac{15}{224}\right)$\tabularnewline
\hline
\end{tabular}
\par\end{center}

Denote quantum dimensions of $L\left(\frac{25}{28},h\right)$, $h=\frac{3}{4},\frac{15}{2},\frac{165}{32},\frac{5}{32},\frac{13}{4},\frac{57}{32}$
by $d_{i},$ $i=1,\cdots,6$ respectively. Then by Proposition \ref{product property of quantum dimension },
we can express quantum dimensions of all the irreducible $L\left(\frac{25}{28},0\right)$-modules
in terms of $1,x,y,d_{i},i=1,\cdots,6$. Direction calculation gives
\begin{align*}
 & \text{glob}L\left(\frac{25}{28},0\right)\\
 & =1+x^{2}+y^{2}+\sum_{i=1}^{6}d_{i}^{2}+\left(x\cdot\sum_{i=1}^{6}d_{i}\right)^{2}+\left(y\cdot\sum_{i=1}^{6}d_{i}\right)^{2}\\
 & =\left(1+\sum_{i=1}^{6}d_{i}^{2}\right)\cdot\left(1+x^{2}+y^{2}\right).
\end{align*}
From Remark \ref{product property of global dimension} we obtain
\[
\text{glob}\left(\mathcal{V}\otimes L\left(\frac{25}{28},0\right)\right)=\left(1+x^{2}+y^{2}\right)^{2}\cdot\left(1+\left(\frac{1+\sqrt{5}}{2}\right)^{2}\right)\cdot\left(1+\sum_{i=1}^{6}d_{i}^{2}\right).
\]
Note that we also have
\[
q\dim_{\mathcal{V}\otimes L\left(\frac{25}{28},0\right)}\mathcal{U}=1+x^{2}+y^{2}.
\]
Since $\mathcal{U}$ is an extension the vertex operator algebra $\mathcal{V}\otimes L\left(\frac{25}{28},0\right)$,
by Theorem \ref{global dim of VOA and subVOA},
\[
\text{glob}\left(\mathcal{V}\otimes L\left(\frac{25}{28},0\right)\right)=\left(q\dim_{\mathcal{V}\otimes L\left(\frac{25}{28},0\right)}\mathcal{U}\right)^{2}\cdot\text{glob \ensuremath{\mathcal{U}}},
\]
which implies
\begin{flalign*}
\text{glob \ensuremath{\mathcal{U}}} & =\frac{\text{glob}\left(\mathcal{V}\otimes L\left(\frac{25}{28},0\right)\right)}{\left(q\dim_{\mathcal{V}\otimes L\left(\frac{25}{28},0\right)}\mathcal{U}\right)^{2}}\\
 & =\frac{\left(1+x^{2}+y^{2}\right)^{2}\cdot\left(1+\left(\frac{1+\sqrt{5}}{2}\right)^{2}\right)\cdot\left(1+\sum_{i=1}^{6}d_{i}^{2}\right)}{\left(1+x^{2}+y^{2}\right)^{2}}\\
 & =\left(1+\left(\frac{1+\sqrt{5}}{2}\right)^{2}\right)\cdot\left(1+\sum_{i=1}^{6}d_{i}^{2}\right).
\end{flalign*}
Now we consider the quantum dimensions of irreducible $\mathcal{U}$-modules
given in Lemma \ref{Real modules}. By Remark \ref{module str.,dual},
$q\dim_{\mathcal{U}}M^{i}=q\dim_{\mathcal{V}\otimes L\left(\frac{25}{28},0\right)}M_{1}^{i}$.
Apply Proposition \ref{product property of quantum dimension }, easy
calculation gives
\begin{center}
\begin{tabular}{|c|c|c|c|c|c|c|c|c|}
\hline
$ $  & $M^{0}$  & $M^{1}$  & $M^{2}$  & $M^{3}$  & $M^{4}$  & $M^{5}$  & $M^{6}$  & $M^{7}$\tabularnewline
\hline
\hline
$q\dim_{\mathcal{U}}$  & $1$  & $d_{1}$  & $d_{5}$  & $d_{2}$  & $d_{3}$  & $d_{4}$  & $d_{6}$  & $\frac{1+\sqrt{5}}{2}$\tabularnewline
\hline
\end{tabular}
\par\end{center}

\vspace{0.2em}

\begin{center}
\begin{tabular}{|c|c|c|c|c|c|c|}
\hline
 & $M^{8}$  & $M^{9}$  & $M^{10}$  & $M^{11}$  & $M^{12}$  & $M^{13}$\tabularnewline
\hline
\hline
$q\dim_{\mathcal{U}}$  & $\frac{1+\sqrt{5}}{2}\cdot d_{1}$  & $\frac{1+\sqrt{5}}{2}\cdot d_{5}$  & $\frac{1+\sqrt{5}}{2}\cdot d_{1}$  & $\frac{1+\sqrt{5}}{2}\cdot d_{2}$  & $\frac{1+\sqrt{5}}{2}\cdot d_{6}$  & $\frac{1+\sqrt{5}}{2}\cdot d_{3}$\tabularnewline
\hline
\end{tabular}
\par\end{center}

From the above table we find
\begin{align*}
 & \sum_{i=0}^{13}\left(q\dim_{\mathcal{U}}M^{i}\right)^{2}\\
 & =1+\sum_{i=1}^{6}d_{i}^{2}+\left(\frac{1+\sqrt{5}}{2}\right)^{2}\left(1+\sum_{i=1}^{6}d_{i}^{2}\right)\\
 & =\left(1+\left(\frac{1+\sqrt{5}}{2}\right)^{2}\right)\left(1+\sum_{i=1}^{6}d_{i}^{2}\right),
\end{align*}
which exactly equals $\text{glob}\mathcal{U}$. Thus these $M^{i},i=0,1,\ \cdots,\ 13$
give all the irreducible modules of $\mathcal{U}$.

Now we obtain the following theorem:

\begin{theorem} $\mathcal{U}$ has exactly 14 inequivalent irreducible
modules, which are listed in Lemma \ref{Real modules}. \end{theorem}

\section{Fusion rules}

In this Section, we shall determine all fusion rules for irreducible
$\mathcal{U}$-modules. We denote by $W^{1}\boxtimes_{\mathcal{U}}W^{2}$
the fusion product of $\mathcal{U}$-modules $W^{1}$ and $W^{2}$,
and $W^{1}\boxtimes_{U^{1}}W^{2}$ the fusion product for $U^{1}$-modules
$W^{1}$ and $W^{2}$.

\begin{theorem} All fusion rules for irreducible $\mathcal{U}$-modules
are given by
\[
\dim_{\mathcal{U}}\left(_{M^{i},M^{j}}^{M^{k}}\right)=\dim_{U^{1}}\left(_{M_{1}^{i},M_{1}^{j}}^{M_{1}^{k}}\right)
\]
where $i,j,k=0,1,\ \cdots,\ 13$.\end{theorem}
\begin{proof}
Since $U^{1}=\mathcal{V}\otimes L\left(\frac{25}{28},0\right)$ is
a rational vertex operator algebra, for irreducible $U^{1}$-modules
$M_{1}^{i}$, $M_{1}^{j}$, $i,j\in\left\{ 0,1,\ \cdots,\ 13\right\} $
we have the fusion product
\begin{equation}
M_{1}^{i}\boxtimes_{U^{1}}M_{1}^{j}=\sum d_{i,j}^{k}W^{k}\label{q-1}
\end{equation}
where $d_{i,j}^{k}=\dim I_{U^{1}}\left(_{M_{1}^{i},M_{1}^{j}}^{W^{k}}\right)$
and $W^{k}$ runs over the set of equivalence classes of irreducible
$U^{1}$-modules. By case by case verification, we find that the fusion
rule $d_{i,j}^{k}=0$ unless $W^{k}\cong M_{1}^{s}$ for some $s\in\left\{ 0,1,\ \cdots,\ 13\right\} $.
Hence $\dim I_{U^{1}}\left(_{M_{1}^{i},M_{1}^{j}}^{M^{k}}\right)=\dim I_{U^{1}}\left(_{M_{1}^{i},M_{1}^{j}}^{M_{1}^{k}}\right)$
and the fusion product in (\ref{q-1}) can be written as
\begin{equation}
M_{1}^{i}\boxtimes_{U^{1}}M_{1}^{j}=\sum_{k=0}^{13}N_{i,j}^{k}M_{1}^{k}\label{eq:4.1}
\end{equation}
where $N_{i,j}^{k}=\dim I_{U^{1}}\left(_{M_{1}^{i},M_{1}^{j}}^{M_{1}^{k}}\right)$.
Since $\mathcal{U}$ is a rational vertex operator algebra, we have
\begin{equation}
M^{i}\boxtimes_{\mathcal{U}}M^{j}=\sum_{k=0}^{13}P_{i,j}^{k}M^{k}\label{4.2}
\end{equation}
where $P_{i,j}^{k}=\dim I_{\mathcal{U}}\left(_{M^{i},M^{j}}^{M^{k}}\right)$.
Remark \ref{module str.,dual} imply
\[
q\dim_{U^{1}}M_{1}^{i}\boxtimes_{U^{1}}M_{1}^{j}=q\dim_{U^{1}}M_{1}^{i}\cdot q\dim_{U^{1}}M_{1}^{j}=q\dim_{\mathcal{U}}M^{i}\cdot q\dim_{\mathcal{U}}M^{j}=q\dim_{\mathcal{U}}M^{i}\boxtimes_{\mathcal{U}}M^{j}.
\]

Then it follows from identities (\ref{eq:4.1}) and (\ref{4.2}) that
\[
\sum_{k=0}^{13}P_{i,j}^{k}q\dim_{\mathcal{U}}M^{k}=\sum_{k=0}^{13}N_{i,j}^{k}q\dim_{U^{1}}M_{1}^{k}.
\]
Note that by Theorem \ref{restriction of fusion rules}, $N_{ij}^{k}\ge P_{ij}^{k}$.
The above equation implies $N_{i,j}^{k}=P_{i,j}^{k}$ and hence the
theorem is proved.
\end{proof}

\vskip10pt {\footnotesize{}{ }\textbf{\footnotesize{}C. Dong}{\footnotesize{}:
Department of Mathematics, University of California Santa Cruz, CA 95064 USA; }\texttt{dong@ucsc.edu}{\footnotesize \par}

\textbf{\footnotesize{}X. Jiao}{\footnotesize{}: Department of Mathematics, East China Normal University, Shanghai 200241, CHINA; }\texttt{xyjiao@math.ecnu.edu.cn}{\footnotesize \par}

\textbf{\footnotesize{}N. Yu}{\footnotesize{}: School of Mathematical
Sciences, Xiamen University, Fujian, 361005, China;} \texttt{
ninayu@xmu.edu.cn}{\footnotesize \par}

{\footnotesize{}}}{\footnotesize \par}

\end{document}